\documentclass[fleqn]{mat01}
\usepackage{times,mathtimy,amssymb,epsfig}
\begin{document}

\setcounter{page}{395}
\firstpage{395}

\makeatletter
\def\artpath#1{\def\@artpath{#1}}
\makeatother \artpath{C:/mathsci-arxiv/nov2003}

\font\xxxx=tir at 9.4pt

\newtheorem{theore}{Theorem}
\renewcommand\thetheore{\arabic{section}.\arabic{theore}}
\newtheorem{propo}[theore]{\rm PROPOSITION}
\newtheorem{lem}[theore]{Lemma}

\newtheorem{theor}[theore]{\bf Theorem}

\def\thoe{\trivlist\item[\hskip\labelsep{{\bf Theorem}}]}
\newtheorem{theo}{\bf Theorem}
\renewcommand\thetheo{\arabic{theo}}
\newtheorem{rem}[theo]{Remark}

\renewcommand{\theequation}{\thesection\arabic{equation}}

\font\zz=msam10 at 10pt
\def\cd{\mbox{\zz{\char'245}}}

\title{Stability estimates for \textbf{h}-\textbf{p}
spectral element methods\\ for general elliptic problems on curvilinear
domains}

\markboth{Pravir Dutt and Satyendra Tomar}{Stability estimates for
\textit{h}-\textit{p} spectral element methods}

\author{PRAVIR DUTT$^{*}$ and SATYENDRA TOMAR$^{\dagger}$}

\address{$^{*}$Department of Mathematics, Indian Institute of
Technology, Kanpur~208~016, India\\
\noindent $^{\dagger}$Department of Applied Mathematics, University of
Twente, P.O.~Box~217, 7500 AE, Enschede, The Netherlands\\
\noindent Email: pravir@iitk.ac.in}

\volume{113}

\mon{November}

\parts{4}


\begin{abstract}
In this paper we show that the h-p spectral element method
developed in \cite{pdstrk1,tomar-01,tomar-dutt-kumar-02} applies to
elliptic problems in curvilinear polygons with mixed Neumann and
Dirichlet boundary conditions provided that the Babuska--Brezzi inf--sup
conditions are satisfied. We establish basic stability estimates for a
non-conforming h-p spectral element method which allows for
simultaneous mesh refinement and variable polynomial degree. The
spectral element functions are non-conforming if the boundary conditions
are Dirichlet. For problems with mixed boundary conditions they are
continuous only at the vertices of the elements. We obtain a stability
estimate when the spectral element functions vanish at the vertices of
the elements, which is needed for parallelizing the numerical scheme.
Finally, we indicate how the mesh refinement strategy and choice of
polynomial degree depends on the regularity of the coefficients of the
differential operator, smoothness of the sides of the polygon and the
regularity of the data to obtain the maximum accuracy achievable.
\end{abstract}

\maketitle

\section{\label{sec:Gen-Intro}Introduction}

In this paper we generalize all the results we have obtained in
\cite{pdstrk1} and seek a numerical solution to an elliptic boundary
value problem where the differential operator satisfies the
Babuska--Brezzi inf--sup conditions. We solve the boundary value problem
on a curvilinear polygon whose sides are piecewise analytic (smooth) and
we assume the boundary conditions are of mixed Neumann and Dirichlet
type as in \cite{babguo1,babguo2,karsher}.

We now briefly describe the contents of this paper. In \S2 we discuss
function spaces and obtain differentiability estimates for the solution
in modified polar coordinates in a sectoral neighbourhood of the
vertices. Here we examine two cases viz. when the coefficients of the
differential operator, sides of the polygon and the data are analytic
and when they have finite regularity.

In \S3 we obtain a stability theorem for a non-conforming spectral
element representation of the solution for problems with mixed boundary
conditions. We let the spectral element functions to be polynomials of
variable degree, where the degree of all these polynomials is bounded by
$W$, and let $M$ denote the number of elements or layers in a sectoral
neighbourhood of each of the vertices in the radial direction as shown
in figure~\ref{Fcursec}. We then define a quadratic form
$\mathcal{V}^{M,W}$ which measures the sum of squares of a weighted
squared norm of the partial differential equation and fractional Sobolev
norms of the boundary conditions and a term which measures the jumps in
the function and its derivatives at inter-element boundaries in
appropriate Sobolev norms. In each of the sectoral neighbourhoods of the
corners we use modified polar coordinates and a global coordinate system
in the remaining part of the domain. We prove that the sum of the
squares of the $H^{2}$ norms of the spectral element functions is
bounded by the quadratic form $\mathcal{V}^{M,W}$ multiplied by a factor
which grows logarithmically in $W$ for problems with Dirichlet boundary
conditions. For problems with mixed boundary conditions this factor can
grow as $M^{4}$, provided $W$ is not too large, and thus the method
displays algebraic instability.

We choose as our approximate solution the unique spectral element
function which minimizes a functional $r^{M,W}$ closely related to the
quadratic form $\mathcal{V}^{M,W}$ as defined in
\cite{pdstrk1,tomar-01,tomar-dutt-kumar-02}. In case the solution is
analytic, we choose $M$ proportional to $W$, and show that $r^{M,W}$
decays exponentially in $M$. Now the error is bounded by $r^{M,W}$
multiplied by a factor which grows at most algebraically in $M$. Hence
the order of convergence remains exponential. If the solution has finite
regularity then we choose $M$ proportional to $\ln W$ and show that
$r^{M,W}$ decays algebraically in $W$. Now the error is bounded by
$r^{M,W}$ multiplied by a factor which grows polylogarithmically in $W$
and hence the error decays algebraically in $W$.

We now come to the aspect of parallelization of the numerical scheme.
For problems with Dirichlet boundary conditions the spectral element
functions are non-conforming and we can use the stability theorem to
parallelize the scheme in an optimal manner. It should be noted that the
method is assymptotically faster then the h-p finite element
method. For problems with mixed boundary conditions we cannot use this
stability theorem to parallelize our method since the factor in the
stability estimate can grow as $M^{4}$. To get around this problem we
make the spectral element functions continuous at the vertices of the
elements only. We then prove a stability theorem for mixed problems when
the spectral element functions vanish at the vertices of their elements.
The values of the spectral element functions at the vertices of their
elements constitute the set of common boundary values we have to solve
for. It should be noted that the cardinality of the set of common
boundary values is much smaller than for finite element methods where
the functions have to be continuous along the edges of the elements.
Since the cardinality of the set of common boundary values is small we
can construct an accurate approximation to the Schur complement matrix
from its definition. As a result the method is faster than the standard
h-p finite element method \cite{tomar-01}.

\section{\label{sec:Gen-Func-spaces-Diff-est}Function spaces and
differentiability estimates}

Let $\Omega $ be a curvilinear polygon with vertices $A_{1},A_{2},
\ldots,A_{p}$ and corresponding sides $\Gamma _{1},\Gamma _{2},\ldots,\Gamma
_{p}$ where $\Gamma _{i}$ joins the points $A_{i-1}$ and $A_{i}$ . We
shall assume that the sides $\overline{\Gamma }_{i}$ are analytic
(smooth) arcs, i.e.
\setcounter{equation}{0}
\begin{equation*}
\overline{\Gamma }_{i}= \{ (\varphi _{i}(\xi ),\
\psi_{i}(\xi ))|\xi \in
\overline{I}=[-1,1].\}
\end{equation*}
with $\varphi_{i}(\xi )$ and $\psi _{i}(\xi )$
being analytic (smooth) functions on $\overline{I}$ and $|\varphi
_{i}^{\prime }(\xi )|^{2}+|\psi _{i}^{\prime
}(\xi )|^{2}\geq \alpha >0$. By $\Gamma _{i}$ we mean
the open arc, i.e. the image of $I=(-1,1)$. Let the angle
subtended at $A_{j}$ be $\omega _{j}.$ We shall denote the boundary
$\partial \Omega $ of $\Omega $ by $\Gamma$. Further let $\Gamma
=\Gamma ^{[0]}\bigcup \Gamma ^{[1]},$ $\Gamma
^{[0]}=\bigcup _{i\in \mathcal{D}}\overline{\Gamma }_{i},$
$\Gamma ^{[1]}=\bigcup _{i\in \mathcal{N}}\overline{\Gamma
}_{i}$ where $\mathcal{D}$ is a subset of the set $\{ i\mid
i=1,\ldots ,p\} $ and $\mathcal{N}=\{ i\mid i=1,\ldots
,p\} \setminus \mathcal{D}.$ Let $x$ denote the vector
$x=(x_{1},x_{2})$.

Let $\mathfrak{L}$ be a strongly elliptic operator
\begin{equation}
\mathfrak{L}(u) =-\sum _{r,s=1}^{2} (a_{r,s}(x)u_{x_{s}})_{x_{r}}+\sum
_{r=1}^{2}b_{r} (x)u_{x_{r}} + c(x)u, \label{eq:Eelloper}
\end{equation}
where $a_{s,r}(x) = a_{r,s}(x),$ $b_{r}(x),$ $c_{r}(x)$ are analytic
(smooth) functions on $\overline{\Omega }$ and for any $(\xi
_{1},\xi _{2})\in \Bbb {R}$ and any $x\in \overline{\Omega }$,
\begin{equation}
\sum _{r,s=1}^{2}a_{r,s}\xi _{r}\xi _{s}\geq \mu _{0} (\xi
_{1}^{2}+\xi _{2}^{2})\label{eq:Eellcond}
\end{equation}
with $\mu _{0}>0$. Moreover let the bilinear form induced by the
operator $\mathfrak{L}$ satisfy the inf--sup conditions.

In this paper we shall consider the boundary value problem
\begin{alignat}{2}
&\mathfrak{L}u = f & &\textrm{on }\Omega,\nonumber \\
&u = g^{[0]} & &\textrm{on }\Gamma
^{[0]},\nonumber \\
&\left(\frac{\partial u}{\partial N} \right)_{A} = g^{[1]}\qquad
& &\textrm{on }\Gamma ^{[1]},\label{eq:Eellbvp}
\end{alignat}
where $\left(\frac{\partial u}{\partial N}\right)_{A}$ denotes the usual
conormal derivative which we shall now define. Let $A$ denote the
$2\times 2$ matrix whose entries are given by
\begin{equation*}
A_{r,s}(x) = a_{r,s} (x)
\end{equation*}
for $r,s=1,2.$ Let $N=(N_{1},N_{2})$ denote the outward
normal to the curve $\Gamma _{i}$ for $i\in \mathcal{N}.$ Then
$\left(\frac{\partial u}{\partial N}\right)_{A}$ is defined as follows:
\begin{equation}
\left(\frac{\partial u}{\partial
N}\right)_{A}(x)=\sum _{r,s=1}^{2}N_{r}a_{r,s}\frac{\partial
u}{\partial x_{s}}.
\end{equation}
We shall assume that the given data $f$ is analytic (smooth) on
$\overline{\Omega }$ and $g^{[l]}$ is analytic (smooth) on every closed
arc $\overline{\Gamma }_{i}$ and $g^{[0]}$ is continuous on
$\Gamma ^{[0]}$.

We need to state our regularity estimates in terms of local variables
which are defined on a geometrical mesh imposed on $\Omega $ as in
\S5 of \cite{babguo2}. We first divide $\Omega $ into subdomains.
Thus we divide $\Omega $ into $p$ subdomains $S^{1},\ldots ,S^{p},$
where $S^{i}$ denotes a domain which contains the vertex $A^{i}$
and no other, and on each $S^{i}$ we define a geometrical mesh. Let
$\mathfrak{S}^{k}= \{ \Omega _{i,j}^{k},j=1,\ldots ,J_{k},i=1,\ldots ,I_{k,j}\} $
be a partition of $S^{k}$ and let $\mathfrak{S}=\bigcup _{k=1}^{p}\mathfrak{S}^{k}.$
Here $J_{k}=M+O(1)$ and $I_{k,j}\leq I$ for all $k$ and $j$, where
$I$ is a constant. As has been stated earlier $M$ denotes the number
of elements or layers in a sectoral neighbourhood of each of the vertices
in the radial direction.

We now put some restrictions on $\mathfrak{S}.$ Let $(r_{k},\theta
_{k})$ denote polar coordinates with center at $A_{k}.$ Let $\tau
_{k}=\ln r_{k}.$ We choose $\rho $ so that the curvilinear sector
$\Omega ^{k}$ with sides $\Gamma _{k}$ and $\Gamma _{k+1},$ center at
$A_{k}$ and radius $\rho $ satisfies
\begin{equation*}
\Omega ^{k}\subseteq \bigcup _{\Omega _{i,j}^{k}\in
\mathfrak{S}_{k}}\overline{\Omega }_{i,j}^{k}.
\end{equation*}
$\Omega ^{k}$ may be represented as
\begin{equation}
\Omega ^{k}=\{ (x,y)\in \Omega :0<r_{k}<\rho \}.
\end{equation}

The geometrical mesh we have imposed on $\Omega $ is as shown in
figure~\ref{Fcursec}.

\begin{fig}
\hskip 4pc{\epsfxsize=5cm\epsfbox{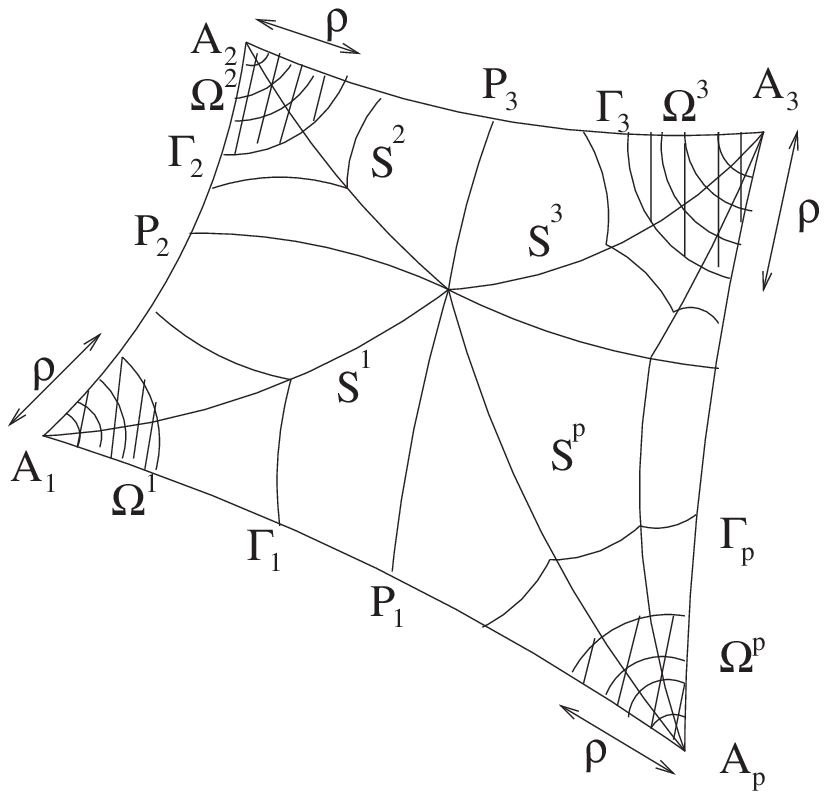}}\vspace{-.5pc}
\caption{\label{Fcursec}Geometric mesh with $M$ layers in the radial
direction in the curvilinear\break}\vspace{-1pc}
\hskip 4pc{\xxxx domain.}\vspace{.5pc}
\end{fig}

\begin{fig}[b]
\hskip 4pc{\epsfxsize=10cm\epsfbox{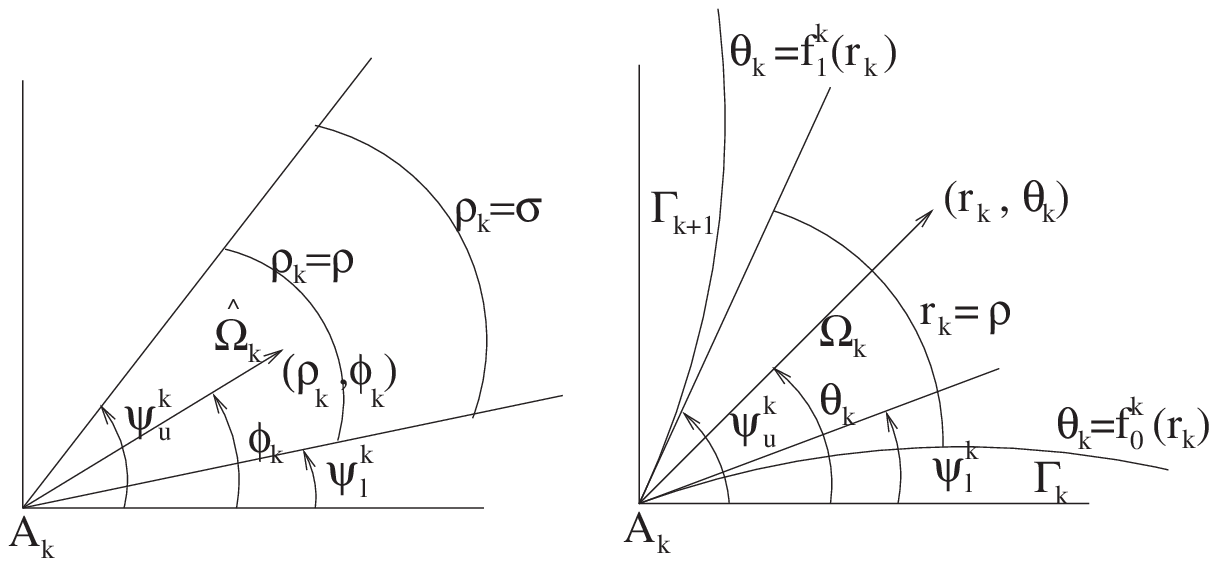}}\vspace{-.5pc}
\caption{\label{F2ksect}Curvilinear sectors.}\vspace{.5pc}
\end{fig}

Let $\gamma _{i,j,l}^{k},$ $1\leq l\leq 4$ be the side of the
quadrilateral $\Omega _{i,j}^{k}\in \mathfrak{S}.$ Then we assume that
\begin{subequations}
\begin{align}
&\gamma _{i,j,l}^{k}: \left\{ \begin{array}{ccc}
 x & = & h_{i,j}^{k}\varphi _{i,j,l}^{k}(\xi ),\\[.5pc]
 y & = & h_{i,j}^{k}\psi _{i,j,l}^{k}(\xi ),\end{array}
\right.\ \ 0\leq \xi \leq 1,l=1,3\label{eqn2.6a}\\[.2pc]
&\gamma _{i,j,l}^{k}: \left\{ \begin{array}{ccc}
 x & = & h_{i,j}^{k}\varphi _{i,j,l}^{k}(\eta ),\\[.5pc]
 y & = & h_{i,j}^{k}\psi _{i,j,l}^{k}(\eta ),\end{array}
\right.\ \ 0\leq \eta \leq 1,l=2,4\label{eqn2.6b}
\end{align}
\end{subequations}
and that for some $C\geq 1$ and $L\geq 1$ independent of $i,j,k$ and
$l$
\begin{equation}
\left|\frac{{\rm d}^{t}}{{\rm d}s^{t}}\varphi
_{i,j,l}^{k}(s)\right|,\left|\frac{{\rm d}^{t}}{{\rm d}s^{t}}\psi
_{i,j,l}^{k}(s)\right|\leq CL^{t}t!,t=1,2,\ldots .\label{eqn2.7}
\end{equation}
We shall also examine the case when they are smooth. Some of the
elements may be triangles too \cite{schwab}. We shall place further
restrictions on the geometric mesh we impose on $\Omega^{k}$ later.

Let $(r_{k},\theta _{k})$ be polar coordinates with center at
$A_{k}.$ Then $\Omega ^{k}$ is the open set bounded by the curvilinear
arcs $\Gamma _{k},$ $\Gamma _{k+1}$ and a portion of the circle
$r_{k}=\rho .$ We subdivide $\Omega ^{k}$ into curvilinear rectangles by
drawing $M$ circular arcs $r_{k}=\sigma _{j}^{k}=\rho \mu
_{k}^{M+1-j},j=2,\ldots ,M+1,$ where $\mu _{k}<1$ and $I_{k}-1$ analytic
curves $C_{2},\ldots ,C_{I_{k}}$ whose exact form we shall prescribe in
what follows. We define $\sigma _{1}^{k}=0.$ Thus $I_{k,j}=I_{k}$ for
$j\leq M;$ in fact, we shall let $I_{k,j}=I_{k}$ for $j\leq M+1.$
Moreover $I_{k,j}\leq I$ for all $k,j$ where $I$ is a fixed constant.
Let
\begin{equation*}
\Gamma _{k+j}= \{(r_{k},\theta _{k})|\theta
_{k}=f_{j}^{k}(r_{k}),0<r_{k}<\rho \},
\end{equation*}
$j=0,1$ in a neighbourhood $A_{k}$ of $\Omega ^{k}$. Then the mapping
\begin{equation}
r_{k}=\rho_{k},\quad \theta _{k}=\frac{1}{(\psi _{u}^{k}-\psi
_{l}^{k})}[(\phi _{k}-\psi
_{l}^{k})f_{1}^{k} (\rho _{k})-(\phi _{k}-\psi
_{u}^{k})f_{0}^{k} (\rho _{k})],\label{eq:Egamdef}
\end{equation}
where $f_{j}^{k}$ is analytic in $r_{k}$ for $j=0,1$, maps locally
the cone
\begin{equation*}
\{ (\rho _{k},\phi _{k}) : 0 < \rho _{k}<\sigma ,\psi
_{l}^{k}<\phi _{k}<\psi _{u}^{k}\}
\end{equation*}
onto a set containing $\Omega ^{k}$ as in \S3 of \cite{babguo2}. The
functions $f_{j}^{k}$ satisfy $f_{0}^{k}(0)=\psi _{l}^{k},$
$f_{1}^{k}(0)=\psi _{u}^{k}$ and
$(f_{j}^{k})^{'}(0)=0$ for $j=0,1.$ It is easy to
see that the mapping defined in (\ref{eq:Egamdef}) has two bounded
derivatives in a neighbourhood of the origin which contains the closure
of the open set
\begin{equation*}
\widehat{\Omega }^{k}= \{(\rho _{k},\phi _{k}):0<\rho
_{k}<\rho ,\psi _{l}^{k}<\phi _{k}<\psi _{u}^{k}\}.
\end{equation*}
We choose the $I_{k-1}$ curves $C_{2},\ldots ,C_{I_{k}}$ as
\begin{equation*}
C_{i}:\phi _{k}(r_{k},\theta _{k})=\psi _{i}^{k}
\end{equation*}
for $i=2,\ldots ,I_{k}.$ Here $\psi _{l}^{k}=\psi _{1}^{k}<\psi
_{2}^{k}<\cdots <\psi _{I_{k}+1}^{k}=\psi _{u}^{k}.$ Let $\Delta \psi
_{i}^{k}=\psi _{i+1}^{k}-\psi _{i}^{k}.$ Then we choose $\{ \psi
_{i}^{k}\} _{i,k}$ so that
\begin{equation}
\max _{i,k}(\Delta \psi _{i}^{k} )<\lambda (\min
_{i,k} (\Delta \psi _{i}^{k}))
\end{equation}
for some constant $\lambda$. We need another set of local variables
$(\tau _{k},\theta _{k})$ in a neighbourhood of $\Omega ^{k}$
where $\tau _{k}=\ln r_{k}.$ In addition we need one final set of
local variables $(\nu _{k},\phi _{k})$ in the cone
\begin{equation*}
\{ (\rho _{k},\phi _{k}):0\leq \rho _{k}\leq \rho ,\psi
_{l}^{k}\leq \phi _{k}\leq \psi _{u}^{k}\},
\end{equation*}
where $\nu _{k}=\ln \rho _{k}.$ Let $S_{\mu }^{k}=\{
(r_{k},\theta _{k}):\: 0\leq r_{k}\leq \mu \} \cap
\Omega .$ Then the image $\widehat{S}_{\mu }^{k}$ in $(\nu
_{k},\phi _{k})$ variables of $S_{\mu }^{k}$ is given by
\begin{equation*}
\widehat{S}_{\mu }^{k}=\{ (\nu _{k},\phi _{k}):\:
-\infty \leq \nu _{k}\leq \ln \mu ,\psi _{l}^{k}\leq \phi ^{k}\leq \psi
_{u}^{k}\}.
\end{equation*}
Now the relationship between the variables $(\tau _{k},\theta _{k})$ and
$(\nu _{k},\phi _{k})$ is given by $(\tau _{k},\theta _{k})=M^{k}(\nu
_{k},\phi _{k})$, viz.
\begin{align}
\tau _{k} &= \nu _{k},\nonumber\\
\theta _{k} &= \frac{1}{(\psi _{u}^{k}-\psi
_{l}^{k})}[(\phi _{k}-\psi
_{l}^{k})f_{1}^{k}({\rm e}^{\nu _{k}})-(\phi _{k}-\psi
_{u}^{k})f_{0}^{k}({\rm e}^{\nu
_{k}})].\label{eq:M-mapping}
\end{align}
Hence it is easy to see that $J_{M^{k}}(\nu _{k},\phi _{k})$,
the Jacobian of the above transformation, satisfies $C_{1}\leq
|J_{M^{k}}(\nu _{k},\phi _{k})|\leq C_{2}$ for all
$(\nu _{k},\phi _{k})\in \widehat{S}_{\mu }^{k},$ for all
$0<\mu \leq \rho .$

We should mention here that it is not necessary to choose the system
of curves we have chosen to impose a geometric mesh on $S_{\mu }^{k}$.
However it is necessary to choose the curve $r_{k}=\rho $ as the
boundary of $\Omega ^{k}$ and no other, as will become apparent in
what follows. Any other additional set of analytic curves which imposes
a geometrical mesh on $S_{\mu }^{k}$ would do equally well. However
the set of curves we have chosen is, in some sense, the most natural
as the image $\widehat{\Omega }_{i,j}^{k}$ of a curvilinear rectangle
$\Omega _{i,j}^{k}$ for $j\geq 2$ in $(\nu _{k},\phi _{k})$
variables is given by a rectangle with straight lines for sides and
for $j=1$ is a semi-infinite strip with straight lines for sides.

We now state the differentiability estimates for the solution $u$
of (\ref{eq:Eellbvp}) which will be needed in this paper.

\begin{propo}\label{4SPdifest}$\left.\right.$\vspace{.5pc}

\noindent Consider the case when the coefficients of the differential
operator are analytic on $\overline{{\Omega }}$ and the sides of the
curvilinear polygon are analytic. Moreover let the geometric mesh
satisfy {\rm (\ref{eqn2.7})}. Let the data $f$ be analytic on
$\overline{{\Omega }}$ and let $g^{[l]}$ be analytic on every closed arc
$\overline{{\Gamma }_{i}}^{[l]}${\rm ,} for $l=0,1${\rm ,} and let
$g^{[0]}$ be continuous on $\Gamma ^{[0]}$. Let $U_{i,j}^{k}(\nu
_{k},\phi _{k})=u(\nu _{k},\phi _{k})$ for $(\nu _{k},\phi _{k})\in
\widehat{\Omega }_{i,j}^{k}$ for $j\leq M$ and $a_{k}=u(A_{k})$. Now
there is an analytic mapping $M_{i,j}^{k}:Sarrow \Omega _{i,j}^{k}$ for
$j>M$ given by $M_{i,j}^{k}(\xi ,\eta )= (X_{i,j}^{k}(\xi ,\eta
),Y_{i,j}^{k}(\xi ,\eta )).$ Here $S$ is the unit square. Let
$U_{i,j}^{k}(\xi ,\eta )=u (X_{i,j}^{k}(\xi ,\eta ),Y_{i,j}^{k}(\xi,
\eta ))$. Then we can show as in {\rm \cite{pdstrk1,tomar-01}} that
\begin{subequations}
\begin{equation}
\Vert U_{i,j}^{k}(\nu _{k},\phi _{k})-a_{k}\Vert
_{m,\widehat{\Omega }_{i,j}^{k}}^{2}\leq (Cm!d^{m}\mu
_{k}^{(1-\beta _{k})(M-j+2)})^{2}
\end{equation}
for $1\leq j\leq M, k=1,\ldots ,p, 1\leq i\leq I_{k}$ and
\begin{equation}
\Vert U_{i,j}^{k}(\xi ,\eta )\Vert _{m,S}^{2}\leq
(Cm!d^{m})^{2}
\end{equation}
\end{subequations}
for $M<j\leq J_{k},$ $1\leq i\leq I_{k,j},$ $1\leq k\leq p.$ Here $C,d$
and $\beta _{k}$ are constants and $0<\beta _{k}<1$ for $1\leq k\leq p$.
\hfill $\cd$
\end{propo}

We next consider the case when the data has finite regularity. To state
the differentiability results in this case we shall need to use the
space $H_{\beta }^{k,l}(\Omega )$ with $k\geq l$ defined in
\cite{babguo1}. We now cite Remark 3 after Theorem 2.1 of
\cite{babguo1}. Let $\overline{{\Gamma }_{j}}$ $\in $$C^{m+2}(\bar{I)}$
for $j=1,\ldots ,p$ and let the coefficients of the differential
operator $\in $$C^{m}(\overline{{\Omega }})$. Let $g^{[0]}$ $\in
$$H_{\beta }^{m+\frac{3}{2},\frac{3}{2}}(\Gamma ^{[0]})$, $g^{[1]}\in
H_{\beta }^{m+\frac{1}{2},\frac{1}{2}}(\Gamma ^{[1]})$ and $f\in
H_{\beta }^{m,0}(\Omega )$. Then there exists a constant
$K_{m}$ such that
\begin{equation}
\Vert u\Vert _{H_{\beta }^{m+2,2}(\Omega )}\leq
K_{m} \left(\Vert f\Vert _{H_{\beta }^{m,0}(\Omega
)}+\sum _{j=0}^{1}\Vert g^{[j]}\Vert _{H_{\beta
}^{m+\frac{3}{2}-j,\frac{3}{2}-j}(\Gamma
^{[j]})} \right).\label{eqn2.12}
\end{equation}

\begin{propo}\label{prop2.2}$\left.\right.$\vspace{.5pc}

\noindent Consider the case when the differential operator and data
satisfy the conditions stated above. We assume moreover that the curves
$\phi _{i,j,l}^{k}$ and $\psi _{i,j,l}^{k}$ defined in {\rm (\ref{eqn2.6a}),
(\ref{eqn2.6b})} satisfy
\begin{equation*}
\Vert \phi _{i,j,l}^{k}\Vert _{m+2,\infty,
\bar{I}},\Vert \psi _{i,j,l}^{k}\Vert _{m+2,\infty,
\bar{I}}\leq E_{m+2}
\end{equation*}
where $E_{m+2}$ is a constant independent of $i,j,k$ and $l$. Let
$U_{i,j}^{k}(\nu _{k},\phi _{k})=u(\nu _{k},\phi
_{k})$ for $(\nu _{k},\phi _{k})\in \hat{\Omega
}_{i,j}^{k}$ for $j\leq M$ and $a_{k}=u(A_{k})$. Now there is
a smooth mapping $M_{i,j}^{k}:Sarrow \Omega _{i,j}^{k}$ for $j>M$
given by $M_{i,j}^{k}(\xi ,\eta )=(X_{i,j}^{k}(\xi
,\eta ),Y_{i,j}^{k}(\xi ,\eta ))$. Here $S$ is
the unit square. Let $U_{i,j}^{k}(\xi ,\eta
)=U(X_{i,j}^{k}(\xi ,\eta ),Y_{i,j}^{k}(\xi
,\eta ))$. Then using {\rm (\ref{eqn2.12})} we can show that
\begin{subequations}
\begin{equation}
\Vert U_{i,j}^{k}(\nu _{k},\phi _{k})-a_{k}\Vert
_{m+2,\hat{\Omega }_{i,j}^{k}}^{2}\leq K_{m+2}(\mu
_{k}^{(1-\beta
_{k})(M-j+2)})^{2}\label{eqn2.13a}
\end{equation}
for $1\leq j\leq M, k=1,\ldots, p, 1\leq i\leq I_{k}$ and
\begin{equation}
\Vert U_{i,j}^{k}(\xi ,\eta )\Vert
_{m+2,S}^{2}\leq K_{m+2}\label{eqn2.13b}
\end{equation}
\end{subequations}
for $M<j\leq J_{k}, 1\leq i\leq I_{k,j}, 1\leq k\leq p$. Here
$K_{m+2}$ denotes a constant.\hfill $\cd$
\end{propo}

\section{\label{sec:Gen-Stab-est}Stability estimates}

\subsection{\it Preliminaries}

Let
\setcounter{equation}{0}
\begin{equation}
\mathfrak{L}u=-\sum
_{r,s=1}^{2}(a_{r,s}(x)u_{x_{s}})_{x_{r}}+\sum
_{r= 1}^{2}b_{r}(x)u_{x_{r}}+c(x)u\label{Estrellop}
\end{equation}
be a strongly elliptic operator which satisfies the inf--sup
conditions. Hence there exists a positive constant $\mu _{0}>0$ such
that
\begin{equation*}
\sum _{r,s=1}^{2}a_{r,s}(x)\xi _{r}\xi _{s}\geq \mu
_{0}(\xi _{1}^{2}+\xi _{2}^{2}),
\end{equation*}
for all $x\in \overline{\Omega }$.

Let $H=H_{0}^{1}(\Omega )$ where $w\in H_{0}^{1}(\Omega )$ if $w\in
H^{1}(\Omega )$ and trace$(w)|_{\Gamma ^{[0]}}=0.$ Consider the
bilinear form $B(u,v)$ defined on $H\times H$ as follows:
\begin{equation}
B(u,v) = \int_{\Omega} \left(\sum
_{r,s=1}^{2}a_{r,s}(x)u_{x_{s}}v_{x_{r}}+\sum
_{r=1}^{2}b_{r}(x)u_{x_{r}}v + cuv \right){\rm d}x.
\end{equation}
Then $B(u,v)$ is a continuous mapping from $H\times
Harrow \Bbb {R}$ and there exists a constant $C_{1}$ such that
\begin{equation}
|B(u,v)|\leq C_{1}\Vert u\Vert _{H^{1}(\Omega )}\Vert v\Vert
_{H^{1}(\Omega)}
\end{equation}
for all $u,v\in H_{0}^{1}(\Omega)$. Moreover we assume that
the inf--sup conditions \cite{schwab}
\begin{subequations}
\begin{equation}
\inf _{0\neq u\in H}\sup _{0\neq v\in
H}\frac{B(u,v)}{\Vert u\Vert _{H^{1}(\Omega
)}\Vert v\Vert _{H^{1}(\Omega )}}\geq
C_{2}>0,\;
\end{equation}
and
\begin{equation}
\sup _{u\in H}B(u,v)>0\quad \mathrm{for}\; \mathrm{every}\;
0\neq v\in H
\end{equation}
\end{subequations}
hold. Then for every continuous linear functional $F(v)$
defined on $H_{0}^{1}(\Omega )$ there exists unique $u_{0}\in
H_{0}^{1}(\Omega )$ such that
$B(u_{0},v)=F(v)$ for all $v\in
H_{0}^{1}(\Omega ).$ Moreover, the {\it a priori} estimate
\begin{equation}
\Vert u_{0}\Vert _{H_{0}^{1}(\Omega )}\leq \frac{1}{C_{2}}\sup _{0\neq
v\in H_{0}^{1}(\Omega )}\frac{|F(v)|}{\Vert v\Vert _{H^{1}(\Omega
)}}\label{Eapriori}
\end{equation}
holds.

Now consider the following mixed boundary value problem
\begin{subequations}
\begin{align}
&\mathfrak{L}u =f\quad \textrm{ in }\Omega,\\
&\overline{\gamma }_{0}u = \left. u \right|_{\Gamma ^{[0]}}
=g^{[0]},
\end{align}
and
\begin{align}
&\overline{\gamma }_{1}u = \left. \left(\frac{\partial u}{\partial
N} \right)_{A} \right|_{\Gamma ^{[1]}} =g^{[1]}.
\end{align}
\end{subequations}
Here the conormal derivative $\overline{\gamma }_{1}u$ is defined as
follows. Let $\Gamma _{i}\subseteq \Gamma ^{[1]}$ and let $T$
and $N$ denote the unit tangent vector and unit outward normal at a
point $P$ on $\Gamma _{i}$ which we traverse in the clockwise direction.
Let $T=(T_{1},T_{2})^{t}$ and
$N=(N_{1},N_{2})^{t}.$ Then
\begin{subequations}
\begin{equation}
\left.\overline{\gamma }_{1}u \right|_{\Gamma _{i}}=
\left.\left(\frac{\partial u}{\partial N}\right)_{A}\right|_{\Gamma
_{i}}=\sum _{r,s=1}^{2}N_{r}a_{r,s}\frac{\partial u}{\partial
x_{s}}=N^{t}A\nabla _{x}u.
\end{equation}
In the same way we define the cotangential derivative
\begin{equation}
\left. \left(\frac{\partial u}{\partial T} \right)_{A} \right|_{\Gamma
_{i}}=\sum _{r,s=1}^{2}T_{r}a_{r,s}\frac{\partial u}{\partial
x_{s}}=T^{t}A\nabla _{x}u,
\end{equation}
and the tangential vector
\begin{equation}
\left. \left(\frac{\partial u}{\partial T} \right)\right|_{\Gamma
_{i}}=T^{t}\nabla _{x}u.
\end{equation}
\end{subequations}
We now consider the spectral elements which are not contained in the
sectoral neighbourhoods of the vertices $\Omega ^{k}$ for $k=1,\ldots
,p$. Now $\Omega _{i,j}^{k}\subseteq \Omega ^{k}$ for $1\leq i\leq
I_{k,j}$ and $1\leq j\leq M$. Let
\begin{equation*}
O^{p+1}=\{ \Omega _{i,j}^{k},1\leq k\leq p,M<j\leq J_{k},1\leq
i\leq I_{k,j}\}.
\end{equation*}
Once more $J_{k}=M+O(1)$. We shall relabel the elements of $O^{p+1}$ and
write
\begin{equation*}
O^{p+1}=\{ \Omega _{l}^{p+1},1\leq l\leq L\} .
\end{equation*}

We shall now introduce some notation so that the reader may proceed
directly to the stability theorem \ref{4STN2allu} and examine the
proof later as it is quite involved.\pagebreak

Consider the domain $\Omega _{l}^{p+1}$. Then there is a mapping
$M_{l}^{p+1}$ from the master square $S=(0,1)\times
(0,1)$ to $\Omega _{l}^{p+1}$. Let $J_{l}^{p+1}(\xi
,\eta )$ denote the Jacobian of the transformation $M_{l}^{p+1}$.
We let
\begin{equation*}
u_{l}^{p+1}(\xi ,\eta )=\sum _{j=0}^{W}\ \sum
_{i=0}^{W}h_{i,j}\xi ^{i}\eta ^{j}.
\end{equation*}

We choose the spectral element functions $\{ u_{i,j}^{k}(\nu
_{k},\phi _{k})\} _{i,j,k}$ for $1\leq i\leq I_{k}$, $1\leq
j\leq M$ and $1\leq k\leq p$ to be polynomials of the form
\begin{equation*}
u_{i,j}^{k}(\nu _{k},\phi _{k})=\sum _{s=0}^{W_{j}}\ \sum
_{r=0}^{W_{j}}a_{r,s}\nu _{k}^{r}\phi _{k}^{s}
\end{equation*}
for $j\ne 1$. Here $1\leq W_{j}\leq W$. If $j=1$ we choose
$u_{i,1}^{k}(\nu _{k},\phi _{k})=g_{k}$ where $g_{k}$ is a
constant for $1\leq i\leq I_{k}$. Let $\pi ^{M,W}$ denote the space of
polynomials $\{ \{ u_{i,j}^{k}(\nu _{k},\phi
_{k})\} _{i,j,k},\{ u_{l}^{p+1}$ $(\xi ,\eta
)\} _{l}\} $.

\setcounter{theo}{0}
\begin{rem} {\rm We shall always choose $M=O(W)$. In case the conditions
of Proposition~{\rm \ref{4SPdifest}} are satisfied so that $u$ is analytic
we choose $W=M$. Once we have obtained the numerical solution we
can define a correction to it so that the corrected solution is conforming
and converges to the actual solution exponentially in $M$ in the
$H^{1}(\Omega )$ norm {\rm \cite{tomar-01,tomar-dutt-kumar-02}}.
Thus the error in the $H^{1}(\Omega )$ norm is bounded
by $C{\rm e}^{-bM}$ where $C$ and $b$ are constants. In case $u\in H_{\beta }^{m+2,2}(\Omega )$
we would choose $M$ proportional to $m\ln W$. Once more we can define
a corrected version of the solution so that it is conforming and converges
to the actual solution in the $H^{1}(\Omega )$ norm and
the error is bounded by $C(\ln W)^{3}W^{-m+1}$. Hence
for the method to converge we must have $m\geq 2$.}
\end{rem}

The stability theorem \ref{4STN2allu} holds provided the coefficients of
the differential operator $\in C^{3}(\bar{\Omega })$ and the
curves $\phi _{i,j,l}^{k},\psi _{i,j,l}^{k}$ defined by (\ref{eqn2.6a}),
(\ref{eqn2.6b}) satisfy
\begin{equation*}
\Vert \phi _{i,j,l}^{k}\Vert _{3,\infty ,\bar{I}},\Vert
\psi _{i,j,l}^{k}\Vert _{3,\infty ,\bar{I}}\leq K_{3},
\end{equation*}
where $K_{3}$ is a constant independent of $i,j,k$ and $l$. In this
paper however we prove Theorem~\ref{4STN2allu} assuming that the
coefficients of the differential operator are analytic on
$\overline{{\Omega }}$ and the curves $\phi _{i,j,l}^{k},\psi
_{i,j,l}^{k}$ defined in (\ref{eqn2.6a}), (\ref{eqn2.6b}) are analytic and
satisfy the condition (\ref{eqn2.7}). Now
\begin{equation*}
\int _{\Omega _{l}^{p+1}}\int |\frak{L}u_{l}^{p+1}(x,y)|^{2}{\rm d}x{\rm d}y = \int
_{S}\int |\frak{L}_{l}^{p+1}u_{l}^{p+1}(\xi ,\eta
)|{\rm d}\xi {\rm d}\eta.
\end{equation*}
Here
\begin{equation*}
\frak {L}_{l}^{p+1}u_{l}^{p+1}(\xi ,\eta )=(\frak
{L}u_{l}^{p+1})(x,y)\sqrt{J_{l}^{p+1}}.
\end{equation*}
Now
\begin{align*}
\frak {L}_{l}^{p+1}w &= A_{l}^{p+1}w_{\xi \xi }+2B_{l}^{p+1}w_{\xi \eta
}+C_{l}^{p+1}w_{\eta \eta }+D_{l}^{p+1}w_{\xi }\\[.3pc]
&\quad\ +E_{l}^{p+1}w_{\eta }+F_{l}^{p+1}w,
\end{align*}
where the coefficients of the differential operator are analytic
(smooth) functions of $\xi $ and $\eta $. Let $\widehat{A}_{l}^{p+1}$ be
the unique polynomial which is the orthogonal projection of
$A_{l}^{p+1}$ into the space of polynomials of degree $W$ in $\xi $ and
$\eta $ with respect to the usual inner product in
$H^{2}\!(S)$. We define
$\widehat{B}_{l}^{p+1},\widehat{C}_{l}^{p+1},\widehat{D}_{l}^{p+1},
\widehat{E}_{l}^{p+1}$ and $\widehat{F}_{l}^{p+1}$ in the same way. We
then define
\begin{align*}
(\frak {L}_{l}^{p+1})^{a}w &= \widehat{A}_{l}^{p+1}w_{\xi \xi
}+2\widehat{B}_{l}^{p+1}w_{\xi \eta }+\widehat{C}_{l}^{p+1}w_{\eta \eta
}+\widehat{D}_{l}^{p+1}w_{\xi }\\[.2pc]
&\quad\ +\widehat{E}_{l}^{p+1}w_{\eta}+\widehat{F}_{l}^{p+1}w.
\end{align*}
Now let $\gamma _{l}$ be a side of the element $\Omega _{m}^{p+1}$ and
let it be the image of the side $\xi =0$ under the mapping
$M_{m}^{p+1}$. Clearly
\begin{equation*}
\frac{\partial u_{m}^{p+1}}{\partial x}=(u_{m}^{p+1})_{\xi
}\xi _{x}+(u_{m}^{p+1})_{\eta }\eta _{x}.
\end{equation*}
We now define
\begin{equation*}
\left. \left(\frac{\partial u_{m}^{p+1}}{\partial
x}\right)^{a} \right|_{\gamma _{l}}=((u_{m}^{p+1})_{\xi
}\widehat{\xi }_{x}+(u_{m}^{p+1})_{\eta }\widehat{\eta
}_{x})(0,\eta ).
\end{equation*}
Here $\widehat{\xi }_{x}(0,\eta )$ and $\widehat{\eta
}_{x}(0,\eta )$ are the unique polynomials which are the
orthogonal projections of $\xi _{x}(0,\eta )$ and $\eta
_{x}(0,\eta )$ into the space of polynomials of degree $W$ in
$\xi $ and $\eta $ with respect to the usual inner product in
$H^{2}(I)$. In the same way we can define
$(\partial u_{m}^{p+1}/\partial y)^{a}$ on $\gamma
_{l}$. Now let $\gamma _{l}$ be a side common to $\Omega _{m}^{p+1}$ and
$\Omega _{n}^{p+1}$ and let it be the image of $\xi =0$ under the
mapping $M_{m}^{p+1}$ and the image of $\xi =1$ under the mapping
$M_{n}^{p+1}$.

Let $[w]$ denote the jump in $w$ across $\gamma _{l}$, where
$w$ is a smooth function on $\overline{\Omega }_{m}^{p+1}$ and
$\overline{\Omega }_{n}^{p+1}$. We now define
\begin{align*}
\left\Vert \left[\left(\frac{\partial u}{\partial
x}\right)^{a} \right] \right\Vert _{1/2,\gamma _{l}}^{2} &=\left\Vert
\left(\frac{\partial u_{m}^{p+1}}{\partial x}\right)^{a}(0,\eta
)- \left(\frac{\partial u_{n}^{p+1}}{\partial
x}\right)^{a}(1,\eta)\right\Vert_{1/2,(0,1)}^{2}\\
\intertext{and}
\left\Vert \left[\left(\frac{\partial u}{\partial
y}\right)^{a}\right]\right\Vert _{1/2,\gamma _{l}}^{2} &= \left\Vert
\left(\frac{\partial u_{m}^{p+1}}{\partial y}\right)^{a}(0,\eta
)-\left(\frac{\partial u_{n}^{p+1}}{\partial
y}\right)^{a}(1,\eta )\right\Vert
_{1/2,(0,1)}^{2}.
\end{align*}
Finally we consider a side $\Gamma _{k}$ of the polygonal domain $\Omega
$ as shown in figure~\ref{Fcursec}. Let $\gamma _{l}$ be a side of
$\Omega _{m}^{p+1}$ such that $\gamma _{l}\subseteq \Gamma _{k}$ and
such that $\gamma _{l}$ is the image of $\xi =0$ under the mapping
$M_{m}^{p+1}$ and which maps the master square $S$ to $\Omega
_{m}^{p+1}$. Then we can define $(\partial u_{m}^{p+1}/\partial T)^{a}$
and $(\partial u_{m}^{p+1}/\partial N)_{A}^{a}$ in the same way. Finally
we define
\begin{align*}
\left\Vert \left(\frac{\partial u}{\partial T}\right)^{a} \right\Vert
_{1/2,\gamma _{l}}^{2} &= \left\Vert \left(\frac{\partial
u_{m}^{p+1}}{\partial T}\right)^{a}(0,\eta ) \right\Vert
_{1/2,(0,1)}^{2}\\
\intertext{and}
\left\Vert \left(\frac{\partial u}{\partial N} \right)_{A}^{a}\right\Vert
_{1/2,\gamma _{l}}^{2} &= \left\Vert \left(\frac{\partial
u_{m}^{p+1}}{\partial N} \right)_{A}^{a}(0,\eta )\right\Vert
_{1/2,(0,1)}^{2}.
\end{align*}
Now consider the sectoral domain $\Omega _{k}$. Let us define the
differential operator
\begin{equation*}
\widetilde{\frak {L}}^{k}w(\tau _{k},\theta _{k}) = {\rm
e}^{2\tau_{k}}\frak {L}w(x,y)
\end{equation*}
as in \cite{pdstrk1}. Then
\begin{equation*}
\widetilde{\frak {L}}^{k}w(\tau _{k},\theta _{k})=\alpha
^{k}w_{\tau _{k}\tau _{k}}+2\beta ^{k}w_{\tau _{k}\theta _{k}}+\gamma
^{k}w_{\theta _{k}\theta _{k}}+\delta ^{k}w_{\tau _{k}}+\epsilon
^{k}w_{\theta _{k}}+\mu ^{k}w,
\end{equation*}

$..$\vspace{-1.5pc}

\noindent where the coefficients of $\widetilde{\frak {L}}^{k}$ are analytic
functions of their arguments. Consider the element $\Omega _{i,j}^{k}$
with $1<j\leq M$. Now the image of $\Omega _{i,j}^{k}$ in $(\nu
_{k},\phi _{k})$ coordinates is the rectangle $\widehat{\Omega
}_{i,j}^{k}$. Clearly
\begin{equation*}
\int _{\Omega _{i,j}^{k}}\int (\widetilde{\frak {L}}^{k}w(\tau
_{k},\theta _{k}))^{2}{\rm d}\tau _{k}{\rm d}\theta _{k}=\int
_{\widehat{\Omega }_{i,j}^{k}}\int (\frak {L}_{i,j}^{k}w(\nu
_{k},\phi _{k}))^{2}{\rm d}\nu _{k}{\rm d}\phi _{k}.
\end{equation*}
Here
\begin{equation*}
\frak {L}_{i,j}^{k}w(\nu _{k},\phi _{k})=\widetilde{\frak
{L}}^{k}w(\tau _{k},\theta _{k})\sqrt{J_{M^{k}}(\nu
_{k},\phi _{k})},
\end{equation*}
where $J_{M^{k}}$ denotes the Jacobian of the transformation $M^{k}$
defined in (\ref{eq:M-mapping}). Once more we can define a differential
operator $(\frak {L}_{i,j}^{k})^{a}$ by replacing the
coefficients of $\frak {L}_{i,j}^{k}$ by polynomials of degree $W$ in
$\nu _{k}$ and $\phi _{k}$ which are exponentially close approximation
to them.

Now the highest order terms of the differential operator
$\widetilde{\frak {L}}^{k}$ are given by $\widetilde{\frak {M}}^{k}$,
where
\begin{equation*}
\widetilde{\frak {M}}^{k}w=\sum _{i,j=1}^{2}\frac{\partial }{\partial
y_{i}} \left(\widetilde{a}_{i,j}^{k}\frac{\partial w}{\partial
y_{j}}\right).
\end{equation*}
Here $y_{1}=\tau _{k}$ and $y_{2}=\theta _{k}$. Let $\widetilde{A}^{k}$
denote the $2\times 2$ matrix such that
$\widetilde{A}_{i,j}^{k}=\widetilde{a}_{i,j}^{k}$. Let $\gamma _{l}$ be
a side of the element $\Omega _{i,j}^{k}$ such that $\gamma
_{l}\subseteq \Gamma _{k}$, where $\Gamma _{k}$ is a side of the polygon
$\Omega $. Let $\widetilde{\gamma }_{l}$ be the image of $\gamma _{l}$
in $(y_{1},y_{2})$ coordinates given by
$y_{1}=y_{1}(\sigma )$, and $y_{2}=y_{2}(\sigma
)$. Let $t$ and $n$ denote the unit tangent and normal vector at a
point $P$ on $\widetilde{\gamma }_{l}$. We now define the conormal
derivative
\begin{equation*}
\left(\frac{\partial w}{\partial
n}\right)_{\widetilde{A}^{k}}=n^{t}\widetilde{A}^{k}\nabla _{y}w.
\end{equation*}
Now the transformation $M^{k}$ defined in (\ref{eq:M-mapping}) maps the
rectangle $\widehat{\Omega }_{i,j}^{k}$ to $\widetilde{\Omega
}_{i,j}^{k}$. Once more we can define $(\partial w/\partial
n)_{\widetilde{A}^{k}}^{a}|_{\widehat{\gamma }_{l}}$ by replacing
the coefficients of the first order differential operator $(\partial
w/\partial n)_{\widetilde{A}^{k}}$ by polynomials of degree $W$ in $\nu
_{k}$ which are exponentially close approximations to them. We can now
define $\Vert (\partial w/\partial n)_{\widetilde{A}^{k}}^{a}
\Vert_{1/2,\widehat{\gamma }_{l}}^{2}$ as we have done before.

The reader can now proceed directly to the stability
theorem~\ref{4STN2allu} stated in \S\ref{sec3.3:Gen-stab-est} and
examine the proof later.

\subsection{\it Technical results}

Consider some $\Omega _{l}^{p+1}\in O^{p+1},$ as shown in
figure~\ref{Frectnortan}. Then $\Omega _{l}^{p+1}$ is a curvilinear
quadrilateral whose sides are analytic arcs and the boundary $\partial
\Omega _{l}^{p+1}$ is traversed in the clockwise direction.

\begin{fig*}[hbt]
\hskip 4pc{\epsfxsize=5cm\epsfbox{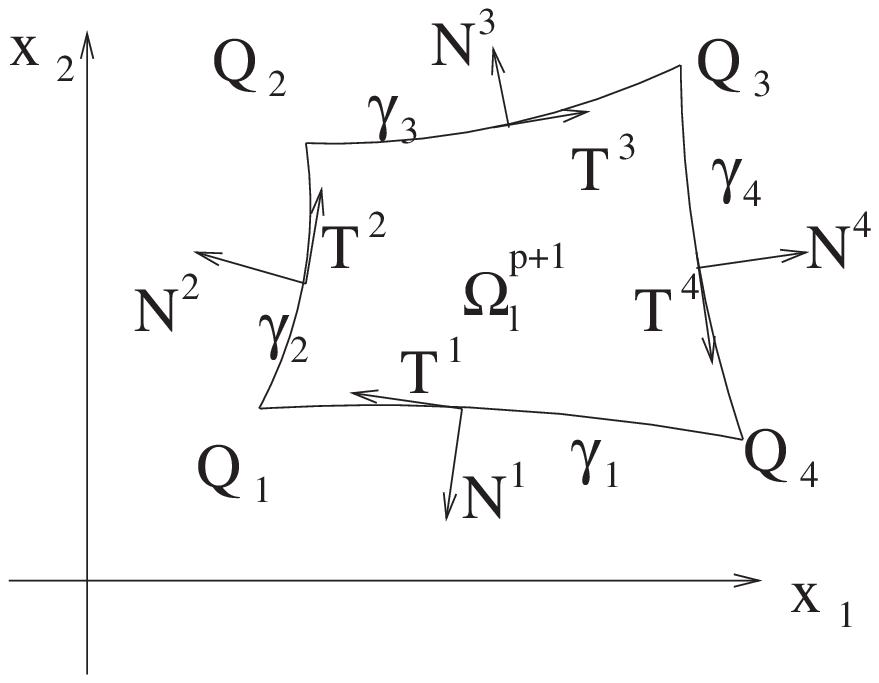}}
\caption{\label{Frectnortan}Element $\Omega _{l}^{p+1}$.}\vspace{.5pc}
\end{fig*}

Let $\gamma $ be a smooth curve and let $N$ and $T$ denote the
unit outward normal and tangent vectors to $\gamma $ at a point $P$
on $\gamma$. Let $s$ be the arc length measured from a point on the
curve in the clockwise direction. Then the second fundamental form
is given by
\begin{equation}
\mathfrak{B}(\xi ,\eta )=-\frac{\partial N}{\partial s}\cdot
T\xi \eta =\frac{\partial T}{\partial s}\cdot N\xi \eta =\kappa \xi
\eta,
\end{equation}
where
\begin{equation*}
\kappa =\pm \frac{{\rm d}T}{{\rm d}s}
\end{equation*}
is the curvature of $\gamma $ at $P$. Clearly
$\mathrm{Trace}(\mathfrak{B})=\kappa$.

Now we need to use Theorem~3.1.1.2 of \cite{grisvard}. Let $v$ be a
smooth vector field defined on $\overline{\Omega }_{l}^{p+1}$ where
$v=(v_{1},v_{2})^{t}.$ Consider the restriction of $v$ to the
boundary $\partial \Omega _{l}^{p+1}.$ Now $\partial \Omega
_{l}^{p+1}=(\bigcup _{i=1}^{4}\gamma _{i})\bigcup
(\bigcup _{i=1}^{4}Q_{i}),$ where $\gamma _{i}$ are the sides
of $\partial \Omega _{l}^{p+1}$ with end points deleted and $Q_{i}$ are
the vertices of $\Omega _{l}^{p+1}.$ We shall denote by $v_{T}$ the
projection of $v$ on the tangent vector $T$ to $\partial \Omega
_{l}^{p+1}$ except at the vertices where this cannot be defined.
Similarly by $v_{N}$ we shall denote the component of $v$ in the
direction of $N.$ Thus we have
\begin{align*}
v_{N} &= v\cdot N\\
\intertext{and}
v_{T} &= v\cdot T.
\end{align*}

\setcounter{theore}{0}
\begin{lem}\label{4SLabsN2u} Let $u\in H^{3}(\Omega_{l}^{p+1})$. Then
\begin{align}
&\frac{\mu _{0}^{2}}{2}\sum _{r,s=1}^{2}\int _{\Omega
_{l}^{p+1}}\left|\frac{\partial ^{2}u}{\partial x_{r}\partial
x_{s}}\right|^{2}{\rm d}x\label{4SELabsN2u}\nonumber\\[.2pc]
&\leq \int _{\Omega _{l}^{p+1}}|\mathfrak{M}u|^{2}{\rm d}x + \sum
_{j=1}^{4}\int _{\gamma _{j}}|\kappa |\left(\left(\frac{\partial u}{\partial
N}\right)_{A}^{2}+\left(\frac{\partial u}{\partial
T}\right)_{A}^{2}\right){\rm d}s\nonumber
\end{align}
\begin{align}
&\quad\ +\frac{512R^{4}}{\mu _{0}^{2}}\sum _{r=1}^{2}\int _{\Omega
_{l}^{p+1}}\left|\frac{\partial u}{\partial x_{r}}\right|^{2}{\rm d}x + 2\sum
_{j=1}^{4}\int _{\gamma _{j}}\left(\frac{\partial u}{\partial
T}\right)_{A}\frac{{\rm d}}{{\rm d}s}\left(\frac{\partial u}{\partial
N}\right)_{A}{\rm d}s\nonumber\\[.2pc]
&\quad\ +\sum _{j=1}^{4}\left\{ \left(\frac{\partial u}{\partial
N^{j+1}}\right)_{A}\left(\frac{\partial u}{\partial
T^{j+1}}\right)_{A}-\left(\frac{\partial u}{\partial
N^{j}}\right)_{A}\left(\frac{\partial u}{\partial
T^{j}}\right)_{A} \right\} (Q_{j}).
\end{align}
\end{lem}
We shall say that a bounded open subset of $\Bbb {R}^{2}$
with Lipschitz boundary $\Gamma $ has a piecewise $C^{2}$ boundary
if $\Gamma =\Gamma _{0}\bigcup \Gamma _{1},$ where
\begin{enumerate}
\renewcommand{\labelenumi}{(\alph{enumi})}
\item $\Gamma _{0}$ has zero measure (for the arc length measure ${\rm d}s$)

\item $\Gamma _{1}$ is open in $\Gamma $ and each point $x\in \Gamma _{1}$
has a $C^{2}$ boundary as defined in 1.2.1.1 of \cite{grisvard}.
Then Theorem~3.1.1.2 of \cite{grisvard} may be stated as follows:
\end{enumerate}

Let $O$ be a bounded open subset of $\Bbb {R}^{2}$ with Lipschitz
boundary $\Gamma$. Assume in addition that $\Gamma $ is piecewise
$C^{2}$. Then for all $v\in (H^{2}(\Omega ))^{2}$
we have
\begin{align}
&\int _{O}|{\rm div}(v)|^{2}{\rm d}x -\int _{O}\sum
_{r,s=1}^{2}\frac{\partial v_{r}}{\partial x_{s}}\frac{\partial
v_{s}}{\partial x_{r}}{\rm d}x\label{Edivv}\nonumber\\[.2pc]
&= \int _{\Gamma _{1}}\left\{ \frac{{\rm d}}{{\rm
d}s}(v_{N}v_{T})-2v_{T}\frac{{\rm d}}{{\rm d}s}v_{N}\right\}
{\rm d}s - \int_{\Gamma_{1}}\{
({\rm tr}\mathfrak{B})v_{N}^{2}+\mathfrak{B}(v_{T},v_{T}
)\}{\rm d}s.
\end{align}
To apply (\ref{Edivv}) we define the vector field
\begin{equation*}
v=A\nabla _{x}u,
\end{equation*}
where $A$ is the matrix
\begin{equation*}
(A)_{r,s}=a_{r,s}.
\end{equation*}
We then observe that
\begin{subequations}
\begin{align}
&\mathfrak{M}u = \sum _{r,s=1}^{2}\frac{\partial }{\partial
x_{r}}\left(a_{r,s}\frac{\partial u}{\partial
x_{s}} \right) = {\rm div}(v),\\[.2pc]
&\left(\frac{\partial u}{\partial N} \right)_{A} = \sum
_{r,s=1}^{2}N_{r}a_{r,s}\frac{\partial u}{\partial
x_{s}}=(\overline{\gamma }_{0}v)\cdot N
\end{align}
and
\begin{equation}
\left(\frac{\partial u}{\partial T}\right)_{A}=\sum
_{r,s=1}^{2}T_{r}a_{r,s}\frac{\partial u}{\partial
x_{s}}=(\overline{\gamma }_{0}v)\cdot T.
\end{equation}
\end{subequations}
\pagebreak

\noindent Hence (\ref{Edivv}) takes the form
\begin{align}
&\int _{\Omega _{l}^{p+1}}|\mathfrak{M}u|^{2}{\rm d}x-\sum
_{r,s=1}^{2}\int _{\Omega _{l}^{p+1}}\frac{\partial v_{r}}{\partial
x_{s}}\frac{\partial v_{s}}{\partial x_{r}}{\rm
d}x\label{Evrxsvsxr}\nonumber\\[.2pc]
&= \sum _{j=1}^{4}\int _{\gamma _{j}}\frac{{\rm d}}{{\rm
d}s}(v_{N}v_{T}){\rm d}s-\sum _{j=1}^{4}2\int
_{\gamma_{j}}\left(\frac{\partial u}{\partial T}\right)_{A}\frac{{\rm
d}}{{\rm d}s}\left(\frac{\partial u}{\partial N}\right)_{A}{\rm
d}s\nonumber\\[.2pc]
&\quad\ -\sum _{j=1}^{4}\int _{\gamma _{j}}\kappa
\left(\left(\frac{\partial u}{\partial
N}\right)_{A}^{2}+\left(\frac{\partial u}{\partial
T}\right)_{A}^{2}\right){\rm d}s.
\end{align}
Now by Lemma 3.1.3.4 of \cite{grisvard} the following inequality holds
for all $u\in H^{2}(\Omega )$:
\begin{equation*}
\mu _{0}^{2}\sum _{r,s=1}^{2}\left|\frac{\partial ^{2}u}{\partial
x_{r}\partial x_{s}}\right|^{2}\leq \sum
_{r,s,k,l=1}^{2}a_{r,k}a_{s,l}\frac{\partial ^{2}u}{\partial
x_{s}\partial x_{k}}\frac{\partial ^{2}u}{\partial x_{r}\partial
x_{l}},
\end{equation*}
a.e. in $\Omega$. Thus it follows that
\begin{equation*}
\mu _{0}^{2}\sum _{r,s=1}^{2}\left|\frac{\partial ^{2}u}{\partial
x_{r}\partial x_{s}}\right|^{2}\leq \sum _{r,s=1}^{2}\frac{\partial
v_{r}}{\partial x_{s}}\frac{\partial v_{s}}{\partial x_{r}}+2\sum
_{r,s,k,l=1}^{2}\left|a_{r,k}\frac{\partial ^{2}u}{\partial
x_{s}\partial x_{k}}\frac{\partial a_{s,l}}{\partial
x_{r}}\frac{\partial u}{\partial x_{l}}\right|,
\end{equation*}
a.e. in $\Omega$. Integrating, we have
\begin{align*}
\mu _{0}^{2}\sum _{r,s=1}^{2}\int\left|\frac{\partial ^{2}u}{\partial
x_{r}\partial x_{s}}\right|^{2}{\rm d}x &\leq \sum _{r,s=1}^{2}\int
\frac{\partial v_{r}}{\partial x_{s}}\frac{\partial v_{s}}{\partial
x_{r}}{\rm d}x\\[.2pc]
&\quad\ + 32R^{2}\int _{\Omega }\sum _{r=1}^{2}\left|\frac{\partial
u}{\partial x_{r}}\right|\sum _{r,s=1}^{2}\left|\frac{\partial
^{2}u}{\partial x_{r}\partial x_{s}}\right|{\rm d}x
\end{align*}
where $R$ is a common bound for all the $C^{1}$ norms of all the
$a_{r,s}$. Hence
\begin{equation}
\frac{\mu _{0}^{2}}{2}\sum _{r,s=1}^{2}\int \left|\frac{\partial
^{2}u}{\partial x_{r}\partial x_{s}}\right|^{2}{\rm d}x\leq \sum
_{r,s=1}^{2}\int \frac{\partial v_{r}}{\partial x_{s}}\frac{\partial
v_{s}}{\partial x_{r}}{\rm d}x + \frac{512R^{4}}{\mu _{0}^{2}}\sum
_{r=1}^{2}\int \left|\frac{\partial u}{\partial
x_{r}}\right|^{2}{\rm d}x.\label{Eu2xrxs}
\end{equation}
Next
\begin{align}
\sum _{j=1}^{4}\int _{\gamma _{j}}\frac{{\rm d}}{{\rm
d}s}(v_{N}v_{T}){\rm d}s &= \sum _{j=1}^{4}\left\{
-\left(\frac{\partial u}{\partial
N^{j+1}}\right)_{A}\left(\frac{\partial u}{\partial
T^{j+1}}\right)_{A} \right.\nonumber\\[.2pc]
&\quad\ \left. + \left(\frac{\partial u}{\partial
N^{j}}\right)_{A}\left(\frac{\partial u}{\partial
T^{j}}\right)_{A}\right\} (Q_{j}).\label{Edervnvt}
\end{align}
Then combining (\ref{Evrxsvsxr})--(\ref{Edervnvt}) we obtain the
result.\hfill $\cd$

\begin{fig*}[b]
\hskip 4pc{\epsfxsize=6cm\epsfbox{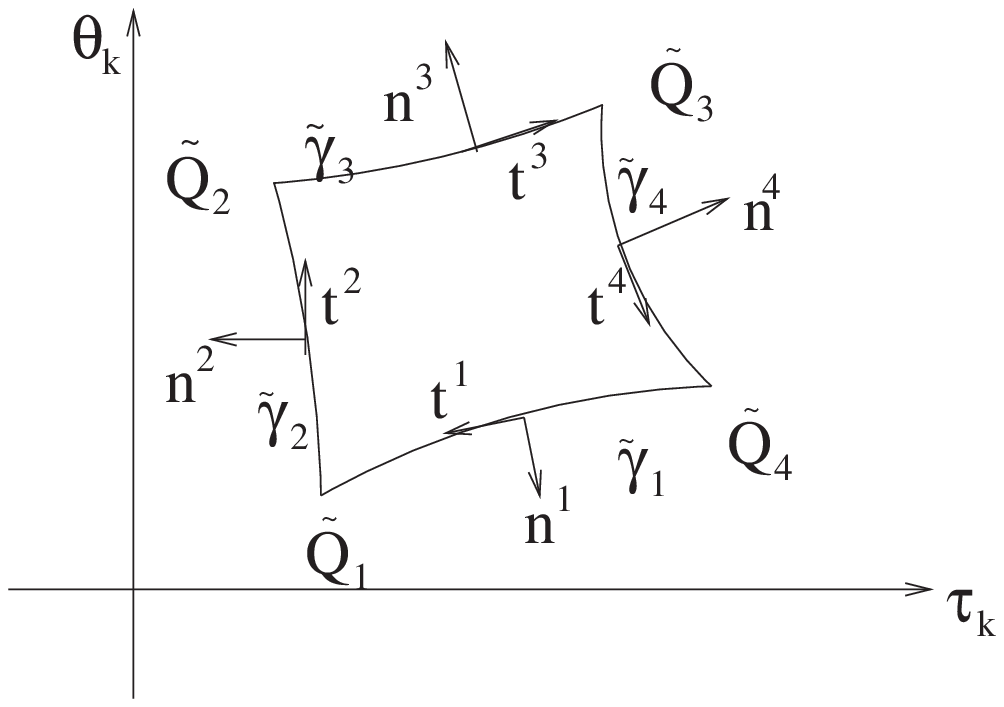}}\vspace{-.5pc}
\caption{\label{Fttrectnortan}Element $\widetilde{\Omega }_{i,j}^{k}$.}\vspace{.5pc}
\end{fig*}

In a neighbourhood of the vertex $A_{k}$ we move to polar coordinates. We
take a curvilinear rectangle $\Omega _{i,j}^{k}$ which comprises part of
the sectoral neighbourhood $\Omega ^{k}$ of the vertex $A_{k}$ and
consider its image $\widetilde{\Omega }_{i,j}^{k}$ in $(\tau
_{k},\theta _{k})$ variables as shown in
figure~\ref{Fttrectnortan}.

As in \cite{pdstrk1} we write the differential operator $\mathfrak{M}$
in modified polar coordinates, where
\begin{equation*}
\mathfrak{M}u=\sum _{r,s=1}^{2}\frac{\partial }{\partial
x_{r}}\left(a_{r,s}\frac{\partial u}{\partial
x_{s}}\right).
\end{equation*}
Now
\begin{align*}
x_{1} &= x_{1}^{k}+{\rm e}^{\tau _{k}}\cos \theta _{k}\\
\intertext{and}
x_{2} &= x_{2}^{k}+{\rm e}^{\tau _{k}}\sin \theta _{k}.
\end{align*}
Here $A_{k}=(x_{1}^{k},x_{2}^{k}).$ We would like to obtain
an estimate for
\begin{equation*}
\int_{\Omega_{i,j}^{k}}r_{k}^{2}|\mathfrak{M}u|^{2}{\rm
d}x=\int _{\widetilde{\Omega
}_{i,j}^{k}}|\widetilde{\mathfrak{M}}^{k}u|^{2}{\rm d}\tau
_{k}{\rm d}\theta _{k}.
\end{equation*}
Let us define the new differential operator
\begin{equation}
\widetilde{\mathfrak{M}}^{k}u={\rm e}^{2\tau _{k}}\sum
_{r,s=1}^{2}\frac{\partial }{\partial x_{r}}\left(a_{r,s}\frac{\partial
u}{\partial x_{s}}\right) = \sum _{r,s=1}^{2}\frac{\partial }{\partial
y_{r}}\left(\tilde{a}_{r,s}\frac{\partial u}{\partial
y_{s}}\right).\label{Emtilku}
\end{equation}
Here $y_{1}=\tau _{k}$ and $y_{2}=\theta _{k}.$ Let $O^{k}$ denote
the matrix
\begin{subequations}
\begin{equation}
O^{k}= \left[\begin{array}{cc}
 \cos \theta _{k} & -\sin \theta _{k}\\[.2pc]
 \sin \theta _{k} & \cos \theta _{k}
\end{array}\right]\label{Edefok}
\end{equation}
and $\widetilde{A}^{k}$ denote the matrix
\begin{equation*}
\widetilde{A}^{k}= \left[\begin{array}{cc}
\widetilde{a}_{1,1}^{k} & \widetilde{a}_{1,2}^{k}\\[.5pc]
\widetilde{a}_{2,1}^{k} & \widetilde{a}_{2,2}^{k}
\end{array}\right].
\end{equation*}
Then it can be easily shown that
\begin{equation}
\widetilde{A}^{k}=(O^{k})^{t}AO^{k}.
\end{equation}
\end{subequations}
Hence, since $O^{k}$ is an orthogonal matrix, we have that
\begin{equation}
\sum _{r,s=1}^{2}\widetilde{a}_{r,s}^{k}\eta _{r}\eta _{s}\geq \mu
_{0}(\eta _{1}^{2}+\eta _{2}^{2}).
\end{equation}
Moreover the following relations hold:
\begin{subequations}
\begin{align}
&(\widetilde{a}_{1,1}^{k})_{\theta _{k}}  =
2\widetilde{a}_{1,2}^{k}+O({\rm e}^{\tau
_{k}}),\label{Eak11theak12}\\[.2pc]
&(\widetilde{a}_{1,2}^{k})_{\theta _{k}}  =
\widetilde{a}_{2,2}^{k}-\widetilde{a}_{1,1}^{k}+O({\rm e}^{\tau
_{k}}),\\[.2pc]
&(\widetilde{a}_{2,2}^{k})_{\theta _{k}}  =
-2\widetilde{a}_{1,2}^{k} + O({\rm e}^{\tau _{k}}),\\[.2pc]
&(\widetilde{a}_{1,1}^{k})_{\tau
_{k}},(\widetilde{a}_{1,2}^{k})_{\tau _{k}}\ \ \
\textrm{and}\ \ \ (\widetilde{a}_{2,2}^{k})_{\tau
_{k}}=O({\rm e}^{\tau _{k}}),\label{Eallaktau}
\end{align}
\end{subequations}
as $\tau_{k}\rightarrow -\infty.$ Next let $\gamma $ be a curve given
by
\begin{align*}
x_{1} & = x_{1}(s),\\[.2pc]
x_{2} & = x_{2}(s),
\end{align*}
where $s$ is the arc length along the curve $\gamma.$ Then the curvature
$\kappa $ at a point $P$ on the curve is given by
\begin{equation*}
\kappa = \frac{{\rm d}x_{1}}{{\rm d}s}\frac{{\rm d}^{2}x_{2}}{{\rm
d}s^{2}}-\frac{{\rm d}x_{2}}{{\rm d}s}\frac{{\rm d}^{2}x_{1}}{{\rm
d}s^{2}}.
\end{equation*}
Let $\widetilde{\gamma }$ be the image of the curve in
$(y_{1},y_{2})$ coordinate given by
\begin{align*}
y_{1} & = y_{1}(\sigma ),\\[.2pc]
y_{2} & = y_{2}(\sigma ),
\end{align*}
where $\sigma $ is the arc length along the curve $\widetilde{\gamma}$.
Then it is easy to verify that
\begin{equation}
\frac{{\rm d}s}{{\rm d}\sigma } = {\rm e}^{y_{1}}.
\end{equation}
Now we can show that the curvature $\widetilde{\kappa }$ of the curve
$\widetilde{\gamma }$ is given by
\begin{equation*}
\widetilde{\kappa} = \kappa {\rm e}^{y_{1}} + \frac{{\rm d}y_{2}}{{\rm
d}\sigma}.
\end{equation*}
Hence
\begin{equation}
|\widetilde{\kappa }|<|\kappa |{\rm e}^{\tau
_{k}} + 1\leq K,
\end{equation}
where $K$ is a uniform constant, for all the curves $\widetilde{\gamma
}_{s}\subseteq \widetilde{\Omega }^{k}$.

We shall denote by $t$ and $n$ the unit tangent and outward normal
vector at a point $P$ on $\widetilde{\gamma},$ the boundary of
$\widetilde{\Omega }_{i,j}^{k}$ except at its vertices where these
are not defined.

\begin{lem}\label{4SLrabsN2u} Let $u(y)\in
H^{3}(\widetilde{\Omega }_{i,j}^{k}).$ Then
\begin{align}
&\frac{\mu _{0}^{2}}{2}\sum _{r,s=1}^{2}\int _{\widetilde{\Omega
}_{i,j}^{k}}\left|\frac{\partial ^{2}u}{\partial y_{r}\partial
y_{s}}\right|^{2}{\rm d}y\label{4SELrabsN2u}\nonumber\\[.2pc]
&\leq \int _{\widetilde{\Omega }_{i,j}^{k}}|\widetilde{\frak
{M}}^{k}u|^{2}{\rm d}y + 2\sum _{j=1}^{4}\int _{\widetilde{\gamma
}_{j}}\left(\frac{\partial u}{\partial
t}\right)_{\widetilde{A}^{k}}\frac{{\rm d}}{{\rm d}\sigma
}\left(\left(\frac{\partial u}{\partial
n}\right)_{\widetilde{A}^{k}}\right){\rm d}\sigma \nonumber\\[.2pc]
&\quad\ +\sum _{j=1}^{4}\left\{\left(\frac{\partial u}{\partial
t^{j+1}}\right)_{\widetilde{A}^{k}}\left(\frac{\partial u}{\partial
n^{j+1}}\right)_{\widetilde{A}^{k}}-\left(\frac{\partial u}{\partial
t^{j}}\right)_{\widetilde{A}^{k}}\left(\frac{\partial u}{\partial
n^{j}}\right)_{\widetilde{A}^{k}}\right\}
(\widetilde{Q}_{j})\nonumber \\[.2pc]
&\quad\ +\sum _{j=1}^{4}\int _{\widetilde{\gamma
}_{j}}|\widetilde{\kappa }|\left(\left(\frac{\partial
u}{\partial t}\right)_{\widetilde{A}^{k}}^{2}+\left(\frac{\partial
u}{\partial n}\right)_{\widetilde{A}^{k}}^{2}\right){\rm d}\sigma
+\frac{512}{\mu _{0}^{2}}R^{4}\sum _{r=1}^{2}\int _{\widetilde{\Omega
}_{i,j}^{k}}\left|\frac{\partial u}{\partial
y_{r}}\right|^{2}{\rm d}y.
\end{align}
\end{lem}

Now once more we use Theorem 3.1.1.2 of \cite{grisvard}. Clearly
$\widetilde{\Omega }_{i,j}^{k}$ for $j\geq 2$ is a bounded open subset
of $\Bbb {R}^{2}$ with Lipschitz boundary $\widetilde{\Gamma }$ that is
a piecewise $C^{2}.$ Thus $\widetilde{\Gamma }=\big(\bigcup
_{i=1}^{4}\widetilde{\gamma }_{i}\big)$ $\bigcup \big(\bigcup
_{i=1}^{4}\widetilde{Q}_{i}\big)$ where $\widetilde{\gamma }_{i}$ are
the sides of the open rectangle $\widetilde{\Omega }_{i,j}^{k}$ with the
end points removed and $\widetilde{Q}_{i}$ are its vertices.

Now
\begin{equation*}
\int _{\Omega _{i,j}^{k}}r_{k}^{2}|\mathfrak{M}u|^{2}{\rm d}x=\int
_{\widetilde{\Omega }_{i,j}^{k}}{\rm e}^{4\tau
_{k}}|\mathfrak{M}u|^{2}{\rm d}\tau _{k}{\rm d}\theta _{k}=\int
_{\widetilde{\Omega
}_{i,j}^{k}}|\widetilde{\mathfrak{M}}^{k}u|^{2}{\rm d}y.
\end{equation*}
Here
\begin{equation*}
\widetilde{\mathfrak{M}}^{k}u=\sum _{r,s=1}^{2}\frac{\partial }{\partial
y_{r}}\left(\widetilde{a}_{r,s}^{k}\frac{\partial u}{\partial
y_{s}}\right)
\end{equation*}
as defined in (\ref{Emtilku}). Then for all $w\in
(H^{2}(\widetilde{\Omega }_{i,j}^{k}))^{2}$ we have
\begin{align}
&\int _{\widetilde{\Omega
}_{i,j}^{k}}|{\rm div}(w)|^{2}{\rm d}y-\sum _{r,s=1}^{2}\int
_{\widetilde{\Omega }_{i,j}^{k}}\frac{\partial w_{r}}{\partial
y_{s}}\frac{\partial w_{s}}{\partial y_{r}}{\rm d}y\label{Edivwotil}\nonumber\\[.2pc]
&= \sum _{j=1}^{4}\left\{ \int _{\widetilde{\gamma
}_{j}}\frac{{\rm d}}{{\rm d}\sigma }(w_{n}w_{t})-2w_{t}\frac{{\rm d}}{{\rm d}\sigma
}w_{n}\right\} {\rm d}\sigma -\sum_{j=1}^{4}\int _{\widetilde{\gamma
}_{j}}\widetilde{\kappa }(w_{n}^{2}+w_{t}^{2}){\rm d}\sigma.
\end{align}
Here $w_{n}$ and $w_{t}$ are the projections of $w$ on the normal
and tangent vectors $n$ and $t$ respectively. We define
\begin{equation*}
w = \widetilde{A}^{k}\nabla _{y}u.
\end{equation*}
Then
\begin{subequations}
\begin{align}
&\widetilde{\mathfrak{M}}^{k}u = \sum _{r,s=1}^{2}\frac{\partial
}{\partial y_{r}}\left(\widetilde{a}_{r,s}^{k}\frac{\partial u}{\partial
y_{s}}\right)={\rm div}(w),\\[.2pc]
&\left(\frac{\partial u}{\partial n}\right)_{\widetilde{A}^{k}} = \sum
_{r,s=1}^{2}n_{r}\widetilde{a}_{r,s}^{k}\frac{\partial u}{\partial
y_{s}}=w_{n},\\
\intertext{and}
&\left(\frac{\partial u}{\partial t}\right)_{\widetilde{A}^{k}} = \sum
_{r,s=1}^{2}t_{r}\widetilde{a}_{r,s}^{k}\frac{\partial u}{\partial
y_{s}}=w_{t}.
\end{align}
\end{subequations}
So (\ref{Edivwotil}) takes the form\vspace{-.2pc}
\begin{align}
&\int _{\widetilde{\Omega
}_{i,j}^{k}}|\widetilde{\mathfrak{M}}^{k}u|^{2}{\rm d}y - \sum
_{r,s=1}^{2}\frac{\partial w_{r}}{\partial y_{s}}\frac{\partial
w_{s}}{\partial y_{r}}{\rm d}y\label{Emtilkuotil}\nonumber\\[.2pc]
&= -2\sum _{j=1}^{4}\int _{\widetilde{\gamma }_{j}}\left(\frac{\partial
u}{\partial t}\right)_{\widetilde{A}^{k}}\frac{{\rm d}}{{\rm d}\sigma
}\left(\left(\frac{\partial u}{\partial
n}\right)_{\widetilde{A}^{k}}\right){\rm d}\sigma\nonumber\\[.2pc]
&\quad\ -\sum _{j=1}^{4}\int
\widetilde{\kappa }\left(\left(\frac{\partial u}{\partial
t}\right)_{\widetilde{A}^{k}}^{2}+\left(\frac{\partial u}{\partial
n}\right)_{\widetilde{A}^{k}}^{2}\right){\rm d}\sigma \nonumber \\[.2pc]
&\quad\ -\sum _{j=1}^{4}\left\{\left(\frac{\partial u}{\partial
t^{j+1}}\right)_{\widetilde{A}^{k}}\left(\frac{\partial u}{\partial
n^{j+1}}\right)_{\widetilde{A}^{k}}-\left(\frac{\partial u}{\partial
t^{j}}\right)_{\widetilde{A}^{k}}\left(\frac{\partial u}{\partial
n^{j}}\right)_{\widetilde{A}^{k}}\right\}
(\widetilde{Q}_{j}).
\end{align}
Now using Lemma~3.1.3.4 of \cite{grisvard} we obtain
\begin{equation*}
\mu _{0}^{2}\sum _{r,s=1}^{2}\left|\frac{\partial ^{2}u}{\partial
y_{r}\partial y_{s}}\right|^{2}\leq \sum _{i,j=1}^{2}\frac{\partial
w_{r}}{\partial y_{s}}\frac{\partial w_{s}}{\partial y_{r}}+2\sum
_{r,s,t,l=1}^{2}\left|\widetilde{a}_{r,t}^{k}\frac{\partial
^{2}u}{\partial y_{s}\partial y_{t}}\frac{\partial
\widetilde{a}_{s,l}^{k}}{\partial y_{r}}\frac{\partial u}{\partial
y_{l}}\right|
\end{equation*}

$\left.\right.$\vspace{-1.5pc}

\noindent and by (\ref{Eak11theak12})--(\ref{Eallaktau}) there exists a
constant $R$ such that $R$ is a common bound for the $C^{1}$ norms of
all $\tilde{a}_{i,j}^{k}$. Hence
\begin{align}
\frac{\mu _{0}^{2}}{2}\sum _{r,s=1}^{2}\int _{\widetilde{\Omega
}_{i,j}^{k}}\left|\frac{\partial ^{2}u}{\partial y_{r}\partial
y_{s}}\right|^{2}{\rm d}y &\leq \sum _{r,s=1}^{2}\int _{\widetilde{\Omega
}_{i,j}^{k}}\frac{\partial w_{r}}{\partial y_{s}}\frac{\partial
w_{s}}{\partial y_{r}}{\rm d}y\nonumber\\[.2pc]
&\quad\ +\frac{512}{\mu _{0}^{2}}R^{4}\sum
_{r=1}^{2}\int _{\widetilde{\Omega }_{i,j}^{k}}\left|\frac{\partial
u}{\partial y_{r}}\right|^{2}{\rm d}y.\label{EabsN2um4}
\end{align}
Thus combining (\ref{Edivwotil}), (\ref{Emtilkuotil}) and
(\ref{EabsN2um4}) we get the result.\hfill $\cd$\vspace{.5pc}

We now need to write terms such as
\begin{equation*}
2\rho^{2}\int _{\gamma _{j}}\left(\frac{\partial u}{\partial
T}\right)_{A}\frac{{\rm d}}{{\rm d}s}\left(\frac{\partial u}{\partial
N}\right)_{A}{\rm d}s
\end{equation*}
in (\ref{4SELrabsN2u}) where $\gamma _{j}\subseteq B_{\rho }^{k}=\{
(x_{1},x_{2}):\rho _{k}=\rho \} $ in terms of
$(y_{1},y_{2})$ coordinates. Let $\gamma $ be a smooth curve
in $\Omega _{\mu }^{k}=\{
(x_{1},x_{2}):(x_{1},x_{2})\in \Omega \;
\mathrm{and}\; \rho _{k}<\mu \} ,$ where $\rho <\mu ,$ and let $P$
be a point on $\gamma $ such that $P$ in polar coordinates has the
representation $(\rho _{k},\theta _{k})$ with $\rho _{k}=\rho$.

Now\vspace{-.2pc}
\begin{equation}
{\rm e}^{y_{1}}\nabla _{x}u=O^{k}\nabla _{y}u,\label{Elapxlapy}
\end{equation}
where $O^{k}$ is the matrix defined in (\ref{Edefok}), and
\begin{equation}
T = O^{k}t,\quad\quad N = O^{k}n.\label{ErelTtNn}
\end{equation}
Hence\vspace{-.2pc}
\begin{subequations}
\begin{equation}
{\rm e}^{y_{1}}\left(\frac{\partial u}{\partial
T}\right)_{A}(P)=t^{t}(O^{k})^{t}AO^{k}\nabla
_{y}u(\widetilde{P})=t^{t}\widetilde{A}^{k}\nabla
_{y}u(\widetilde{P})=\left(\frac{\partial u}{\partial
t}\right)_{\widetilde{A}^{k}}(\widetilde{P})
\end{equation}
using (\ref{Edefok}), (\ref{Elapxlapy}) and (\ref{ErelTtNn}). Here
$\widetilde{P}$ is the image of the point $P$ in
$(y_{1},y_{2})$ coordinates. Similarly, we have
\begin{equation}
{\rm e}^{y_{1}}\left(\frac{\partial u}{\partial
N}\right)_{A}(P)=\left(\frac{\partial u}{\partial
n}\right)_{\widetilde{A}^{k}}(\widetilde{P}).
\label{EuNAunAtil}
\end{equation}
\end{subequations}\vspace{.01pc}

\setcounter{theore}{0}
\begin{propo}\label{4SPutddsun}$\left.\right.$\vspace{.5pc}

\noindent Thus we can conclude that
\begin{subequations}
\begin{equation}
2\rho ^{2}\int _{\gamma _{j}}\left(\frac{\partial u}{\partial
T}\right)_{A}\frac{{\rm d}}{{\rm d}s}\left(\frac{\partial u}{\partial
N}\right)_{A}{\rm d}s=2\int _{\widetilde{\gamma }_{j}}\left(\frac{\partial
u}{\partial t}\right)_{\widetilde{A}^{k}}\frac{{\rm d}}{{\rm d}\sigma
}\left(\frac{\partial u}{\partial n}\right)_{\widetilde{A}^{k}}{\rm d}\sigma
\end{equation}
and\vspace{-.2pc}
\begin{equation}
\left\{ \rho ^{2}\left(\frac{\partial u}{\partial
T}\right)_{A}\left(\frac{\partial u}{\partial N}\right)_{A}\right\}
(P) = \left\{\left(\frac{\partial u}{\partial
t}\right)_{\widetilde{A}^{k}}\left(\frac{\partial u}{\partial
n}\right)_{\widetilde{A}^{k}}\right\} (\widetilde{P}).
\end{equation}
\end{subequations}

\hfill $\cd$\vspace{-.2pc}
\end{propo}

In the same way we obtain the following results.

\begin{propo}\label{4SPubd3}$\left.\right.$\vspace{.5pc}

\noindent Consider the boundary $\gamma $ common to $\Omega
_{i,M+1}^{k}$ and $\Omega _{i,M}^{k}$. Then the following relations hold
{\rm (}figure~{\rm \ref{F2rectnortan}):}
\begin{subequations}
\begin{align}
\left\{ \rho^{2}\left(\frac{\partial u}{\partial
T^{3}}\right)_{A}\left(\frac{\partial u}{\partial
N^{3}}\right)_{A}\right\} (Q_{1}) &= \left\{
\left(\frac{\partial u}{\partial
t^{3}}\right)_{\widetilde{A}^{k}}\left(\frac{\partial u}{\partial
n^{3}}\right)_{\widetilde{A}^{k}}\right\}
(\widetilde{Q}_{1}),\\[.2pc]
\left\{ \rho^{2}\left(\frac{\partial u}{\partial
T^{2}}\right)_{A}\left(\frac{\partial u}{\partial
N^{2}}\right)_{A}\right\} (Q_{1}) &= \left\{
\left(\frac{\partial u}{\partial
t^{4}}\right)_{\widetilde{A}^{k}}\left(\frac{\partial u}{\partial
n^{4}}\right)_{\widetilde{A}^{k}}\right\}
(\widetilde{Q}_{1}),\\[.2pc]
\left\{ \rho^{2}\left(\frac{\partial u}{\partial
T^{2}}\right)_{A}\left(\frac{\partial u}{\partial
N^{2}}\right)_{A}\right\} (Q_{2}) &= \left\{
\left(\frac{\partial u}{\partial
t^{4}}\right)_{\widetilde{A}^{k}}\left(\frac{\partial u}{\partial
n^{4}}\right)_{\widetilde{A}^{k}}\right\}
(\widetilde{Q}_{2}),
\end{align}
and\vspace{-.2pc}
\begin{equation}
\left\{ \rho ^{2}\left(\frac{\partial u}{\partial
T^{1}}\right)_{A}\left(\frac{\partial u}{\partial
N^{1}}\right)_{A}\right\} (Q_{2}) = \left\{
\left(\frac{\partial u}{\partial
t^{1}}\right)_{\widetilde{A}^{k}}\left(\frac{\partial u}{\partial
n^{1}}\right)_{\widetilde{A}^{k}}\right\}
(\widetilde{Q}_{2}).
\end{equation}
\end{subequations}

\begin{fig*}
\hskip 4pc{\epsfxsize=10cm\epsfbox{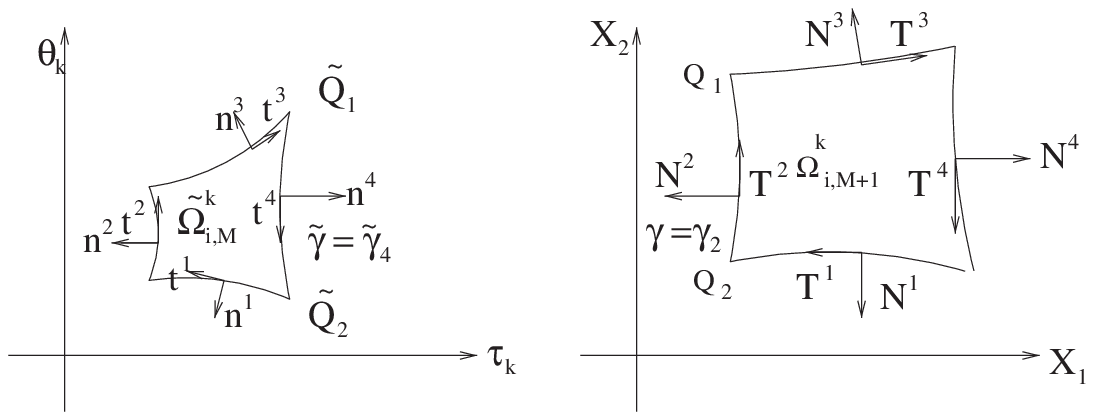}}\vspace{-.5pc}
\caption{\label{F2rectnortan}Elements $\widetilde{\Omega }_{i,M}^{k}$ and
$\Omega _{i,M+1}^{k}$.}\vspace{.5pc}
\end{fig*}

\hfill $\cd$
\end{propo}

Now let $\widetilde{\gamma }_{l}\subseteq \partial \widetilde{\Omega
}_{i,j}^{k}$ for some $j\leq M$ and further suppose $\widetilde{\gamma
}_{l}\subseteq \widetilde{\Gamma }_{j}$ where $j\in \mathcal{D}.$ Let
$n$ and $t$ be the unit outward normal and tangent vectors,
respectively, defined at every point of $\widetilde{\gamma }_{l}$. Then \vspace{-.2pc}
\begin{subequations}
\begin{equation}
\left(\frac{\partial u}{\partial
t}\right)_{\widetilde{A}^{k}}(\sigma
)=\widetilde{g}^{k}(\sigma )\left(\frac{\partial
u}{\partial t}\right)(\sigma )+\widetilde{h}^{k}(\sigma
)\left(\frac{\partial u}{\partial
n}\right)_{\widetilde{A}^{k}}(\sigma ).
\end{equation}

\pagebreak

\noindent Here $\sigma $ is the arc length measured from the point $\widetilde{G}$
(figure~\ref{Fcurpoint}) where
\begin{equation}
\widetilde{g}^{k}(\sigma )=t^{t}\widetilde{A}^{k}t(\sigma
)-\frac{(t^{t}\widetilde{A}^{k}n(\sigma
))^{2}}{n^{t}\widetilde{A}^{k}n(\sigma )},\label{4SEgtilk}
\end{equation}
and
\begin{equation}
\widetilde{h}^{k}(\sigma )=\frac{t^{t}\widetilde{A}^{k}n(\sigma
)}{n^{t}\widetilde{A}^{k}n(\sigma )}.\label{4SEhtilk}
\end{equation}
\end{subequations}
Hence
\begin{align*}
\int _{\widetilde{\gamma }_{l}} \left(\frac{\partial u}{\partial
t}\right)_{\widetilde{A}^{k}}\frac{{\rm d}}{{\rm d}\sigma
}\left(\frac{\partial u}{\partial n}\right)_{\widetilde{A}^{k}}{\rm
d}\sigma &= \int _{\widetilde{\gamma }_{l}}\widetilde{g}^{k}(\sigma
)\frac{\partial u}{\partial t}\frac{{\rm d}}{{\rm d}\sigma
}\left(\frac{\partial u}{\partial n}\right)_{\widetilde{A}^{k}}{\rm
d}\sigma\\[.2pc]
&\quad\ +\int _{\widetilde{\gamma
}_{l}}\frac{\widetilde{h}^{k}(\sigma )}{2}\frac{{\rm d}}{{\rm
d}\sigma }\left(\left(\frac{\partial u}{\partial
n}\right)_{\widetilde{A}^{k}}^{2}\right){\rm d}\sigma.
\end{align*}

And so we can conclude that the following holds.

\begin{propo}\label{4SPudatilk}

\begin{align}
&\int _{\widetilde{\gamma }_{l}}\left(\frac{\partial u}{\partial
t}\right)_{\widetilde{A}^{k}}\frac{{\rm d}}{{\rm d}\sigma
}\left(\frac{\partial u}{\partial n}\right)_{\widetilde{A}^{k}}{\rm
d}\sigma\nonumber\\[.2pc]
&= \int _{\widetilde{\gamma }_{l}}\widetilde{g}^{k}(\sigma
)\frac{\partial u}{\partial t}\frac{{\rm d}}{{\rm d}\sigma
}\left(\frac{\partial u}{\partial n}\right)_{\widetilde{A}^{k}}{\rm
d}\sigma\nonumber\\[.2pc]
&\quad\ -\frac{1}{2}\int _{\widetilde{\gamma }_{l}}\frac{{\rm
d}\widetilde{h}^{k}}{{\rm d}\sigma }\left(\frac{\partial u}{\partial
n}\right)_{\widetilde{A}^{k}}^{2}{\rm d}\sigma
+\left. \frac{\widetilde{h}^{k}(\sigma
)}{2}\left(\frac{\partial u}{\partial
n}\right)_{\widetilde{A}^{k}}^{2}\right|_{\partial \widetilde{\gamma
}_{l}}.
\end{align}
Here $\widetilde{g}^{k}(\sigma )$ and
$\widetilde{h}^{k}(\sigma )$ are defined in {\rm
(\ref{4SEgtilk})} and {\rm (\ref{4SEhtilk})}.\hfill $\cd$
\end{propo}

\begin{fig*}
\hskip 4pc{\epsfxsize=6cm\epsfbox{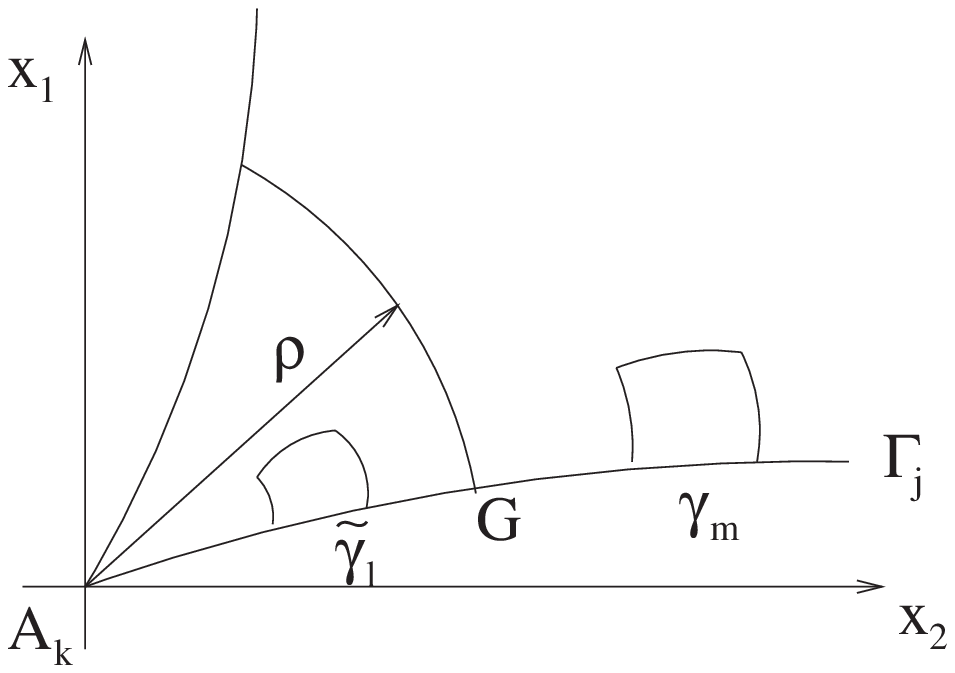}}\vspace{-.5pc}
\caption{\label{Fcurpoint}Arc length measured from the point $G$.}\vspace{.5pc}
\end{fig*}

Next let $\gamma _{m}\subseteq \partial \Omega _{i,j}^{k}$ for some
$j>M$ such that $\gamma _{m}\subseteq \Gamma _{j}$ where $j\in
\mathcal{D}.$ Let $N$ and $T$ be the unit normal and tangent vectors,
respectively, defined at every point of $\gamma _{m}.$ Then
\begin{subequations}
\begin{equation}
\left(\frac{\partial u}{\partial
T}\right)_{A}(s)=g(s)\left(\frac{\partial
u}{\partial T}\right)(s)+h(s)\left(\frac{\partial
u}{\partial N}\right)_{A}(s),
\end{equation}
where $s$ is the arc length measured from the point $G$ as shown in
figure~\ref{Fcurpoint}. Here
\begin{align}
g(s) &= T^{t}AT-\frac{(T^{t}AN)^{2}}{N^{t}AN},\label{4SEgons}\\
\intertext{and}
h(s) &= \frac{T^{t}AN}{N^{t}AN}.\label{4SEhons}
\end{align}
\end{subequations}
So we obtain the following result.

\begin{propo}\label{4SPudongamm}
\begin{align}
&\rho^{2}\int _{\gamma _{m}}\left(\frac{\partial u}{\partial
T}\right)_{A}\frac{{\rm d}}{{\rm d}s}\left(\frac{\partial u}{\partial
N}\right)_{A}{\rm d}s\nonumber\\[.2pc]
&= \rho^{2}\int _{\gamma _{m}}g(s)\frac{\partial u}{\partial
T}\frac{{\rm d}}{{\rm d}s}\left(\frac{\partial u}{\partial
N}\right)_{A}{\rm d}s\nonumber\\[.2pc]
&\quad\ - \frac{\rho ^{2}}{2}\int _{\gamma
_{m}}\frac{{\rm d}h}{{\rm d}s}\left(\frac{\partial u}{\partial
N}\right)_{A}^{2}{\rm d}s + \left.\frac{\rho ^{2}h}{2}\left(\frac{\partial
u}{\partial N}\right)_{A}^{2}\right|_{\partial \gamma _{m}}.
\end{align}

\hfill $\cd$
\end{propo}

Now by (\ref{EuNAunAtil}) we have that
\begin{equation*}
\rho ^{2}\left(\frac{\partial u}{\partial
N}\right)_{A}^{2}(G)=\left(\frac{\partial u}{\partial
n}\right)_{\widetilde{A}^{k}}^{2}(\widetilde{G}).
\end{equation*}
And moreover by (\ref{Edefok}) and (\ref{ErelTtNn})
\begin{subequations}
\begin{align}
g(G) &= \widetilde{g}^{k}(\widetilde{G}),\\
\intertext{and}
h(G) &= \widetilde{h}^{k}(\widetilde{G}).
\end{align}
\end{subequations}

We can now prove the following estimate.

\setcounter{theore}{2}
\begin{lem}\label{4SLabsu012exet}
Let $u_{l}^{p+1}\in H^{3}(\Omega _{l}^{p+1})$. Then
\begin{align}
&\sum_{|\alpha |=2}\int _{S}\int |D_{\xi }^{\alpha
_{1}}D_{\eta }^{\alpha _{2}}u_{l}^{p+1}(\xi ,\eta
)|^{2}{\rm d}\xi {\rm d}\eta\nonumber\\[.2pc]
&\quad\ -C\left(\sum _{|\alpha |\leq
1}\int _{S}\int |D_{\xi }^{\alpha _{1}}D_{\eta }^{\alpha
_{2}}u_{l}^{p+1}|^{2}{\rm d}\xi {\rm d}\eta \right)\nonumber\\[.2pc]
&\leq K\int _{S}\int
|\mathfrak{L}_{l}^{p+1}u_{l}^{p+1}|^{2}{\rm d}\xi {\rm d}\eta +2\rho
^{2}\sum _{r=1}^{4}\int \left(\frac{\partial u_{l}^{p+1}}{\partial
T}\right)_{A}\frac{{\rm d}}{{\rm d}s}\left(\frac{\partial u_{l}^{p+1}}{\partial
N}\right)_{A}{\rm d}s\nonumber \\[.2pc]
&\quad\ +\sum _{r=1}^{4}\rho ^{2}\left\{\left(\frac{\partial
u_{l}^{p+1}}{\partial N^{r+1}}\right)_{A}\left(\frac{\partial
u_{l}^{p+1}}{\partial T^{r+1}}\right)_{A}-\left(\frac{\partial
u_{l}^{p+1}}{\partial N^{r}}\right)_{A}\left(\frac{\partial
u_{l}^{p+1}}{\partial T^{r}}\right)_{A}\right\}
(Q_{r})\nonumber \\[.2pc]
&\quad\ +\sum _{r=1}^{4}\int _{\gamma _{r}}|\kappa |\rho
^{2}\left(\left(\frac{\partial u_{l}^{p+1}}{\partial
N}\right)_{A}^{2}+\left(\frac{\partial u_{l}^{p+1}}{\partial
T}\right)_{A}^{2}\right){\rm d}s.\label{EabsN012upp1l}
\end{align}
Here $S$ is the unit square and $\mathfrak{L}_{l}^{p+1}$ is the
differential operator $\mathfrak{L}$ written in $(\xi ,\eta )$
coordinates. Here $K$ and $C$ are positive constants.
\end{lem}

Recall that
\begin{align}
\mathfrak{L}u &= -\sum
_{r,s=1}^{2}(a_{r,s}(x)u_{x_{s}})_{x_{r}}+\sum
_{r=1}^{2}b_{r}(x)u_{x_{r}}+c(x)u\nonumber\\[.2pc]
&= \mathfrak{M}u+\mathfrak{N}u,
\end{align}
where
\begin{equation*}
\mathfrak{N}u =
\sum_{r=1}^{2}b_{r}(x)u_{x_{r}}+c(x)u.
\end{equation*}
Hence
\begin{equation*}
\rho^{2}\int _{\Omega _{l}^{p+1}}|\mathfrak{M}u|^{2}{\rm d}x\leq
2\rho^{2}\int _{\Omega
_{l}^{p+1}}|\mathfrak{L}u|^{2}{\rm d}x + 2\rho ^{2}\int _{\Omega
_{l}^{p+1}}|\mathfrak{N}u|^{2}{\rm d}x.
\end{equation*}
Using Lemma~\ref{4SLabsN2u} we can conclude that there is
a constant $C$ such that the following estimate holds.
\begin{align}
&\frac{\rho ^{2}\mu _{0}^{2}}{2}\sum _{r,s=1}^{2}\int _{\Omega
_{l}^{p+1}} \left|\frac{\partial ^{2}u_{l}^{p+1}}{\partial x_{r}\partial
x_{s}} \right|^{2}{\rm d}x\nonumber\\[.2pc]
&\quad\ -C\rho ^{2} \left(\sum _{r=1}^{2}\left(\int _{\Omega
_{l}^{p+1}}\left|\frac{\partial u_{l}^{p+1}}{\partial
x_{r}}\right|^{2}{\rm d}x\right)+\int _{\Omega
_{l}^{p+1}}|u_{l}^{p+1}|^{2}{\rm d}x\right)\nonumber
\end{align}
\begin{align}
&\leq 2\rho ^{2}\int _{\Omega
_{l}^{p+1}}|\mathfrak{L}u_{l}^{p+1}|^{2}{\rm d}x+2\rho ^{2}\sum
_{j=1}^{4}\int _{\gamma _{j}}\left(\frac{\partial u_{l}^{p+1}}{\partial
T}\right)_{A}\frac{{\rm d}}{{\rm d}s}\left(\frac{\partial u_{l}^{p+1}}{\partial
N}\right)_{A}{\rm d}s\nonumber \\[.2pc]
&\quad\ +\sum _{j=1}^{4}\rho ^{2}\left\{ \left(\frac{\partial
u_{l}^{p+1}}{\partial N^{j+1}}\right)_{A}\left(\frac{\partial
u_{l}^{p+1}}{\partial T^{j+1}}\right)_{A}-\left(\frac{\partial
u_{l}^{p+1}}{\partial N^{j}}\right)_{A}\left(\frac{\partial
u_{l}^{p+1}}{\partial T^{j}}\right)_{A}\right\}
(Q_{j})\nonumber \\[.2pc]
&\quad\ +\rho ^{2}\sum _{j=1}^{4}\int _{\gamma _{j}}|\kappa
|\left(\left(\frac{\partial u_{l}^{p+1}}{\partial
N}\right)_{A}^{2}+\left(\frac{\partial u_{l}^{p+1}}{\partial
T}\right)_{A}^{2}\right){\rm d}s.\label{EabsN012u}
\end{align}
Writing the above in $(\xi ,\eta )$ coordinates we obtain the
result.\hfill $\cd$\vspace{.5pc}

In the same way we can prove the following estimate.

\begin{lem}\label{4SLabsn2unufi} Let $u_{i,j}^{k}\in H^{3}(\Omega
_{i,j}^{k})$. Then
\begin{align}
&\beta \sum _{|\alpha |=2}\int _{\widehat{\Omega
}_{i,j}^{k}}\int |D_{\nu _{k}}^{\alpha _{1}}D_{\phi _{k}}^{\alpha
_{2}}u_{i,j}^{k}|^{2}{\rm d}\nu _{k}{\rm d}\phi _{k}\nonumber\\[.2pc]
&\quad\ - C\left(\sum _{|\alpha |=1}\int _{\widehat{\Omega
}_{i,j}^{k}}\int |D_{\nu _{k}}^{\alpha _{1}}D_{\phi _{k}}^{\alpha
_{2}}u_{i,j}^{k}|^{2}{\rm d}\nu _{k}{\rm d}\phi _{k}+\int _{\widehat{\Omega
}_{i,j}^{k}}\int |u_{i,j}^{k}|^{2}{\rm e}^{4\nu _{k}}{\rm d}\nu _{k}{\rm d}\phi
_{k}\right)\nonumber\\[.2pc]
&\leq K\int _{\widehat{\Omega }_{i,j}^{k}}\int
|\mathfrak{L}_{i,j}^{k}u_{i,j}^{k}|^{2}{\rm d}\nu _{k}{\rm d}\phi
_{k}+2\sum _{r=1}^{4}\int _{\widetilde{\gamma }_{r}}\left(\frac{\partial
u_{i,j}^{k}}{\partial t}\right)_{\widetilde{A}^{k}}\frac{{\rm d}}{{\rm d}\sigma
}\left(\frac{\partial u_{i,j}^{k}}{\partial
n}\right)_{\widetilde{A}^{k}}{\rm d}\sigma \nonumber \\[.2pc]
&\quad\ +\sum _{r=1}^{4}\left\{ \left(\frac{\partial
u_{i,j}^{k}}{\partial
n^{r+1}}\right)_{\widetilde{A}^{k}}\left(\frac{\partial
u_{i,j}^{k}}{\partial
t^{r+1}}\right)_{\widetilde{A}^{k}}-\left(\frac{\partial
u_{i,j}^{k}}{\partial
n^{r}}\right)_{\widetilde{A}^{k}}\left(\frac{\partial
u_{i,j}^{k}}{\partial t^{r}}\right)_{\widetilde{A}^{k}}\right\}
(\widetilde{Q}_{r})\nonumber \\[.2pc]
&\quad\ +\sum _{r=1}^{4}\int _{\widetilde{\gamma
}_{r}}|\widetilde{\kappa }|\left(\left(\frac{\partial
u_{i,j}^{k}}{\partial
n}\right)_{\widetilde{A}^{k}}^{2}+\left(\frac{\partial
u_{i,j}^{k}}{\partial
t}\right)_{\widetilde{A}^{k}}^{2}\right){\rm d}\sigma.\label{EabsN2unufi}
\end{align}
Here $\widehat{\Omega }_{i,j}^{k}=(\psi _{i}^{k},\psi
_{i+1}^{k})\times (\alpha _{j}^{k},\alpha _{j+1}^{k})$
and $\beta, C$ and $K$ are positive constants.
\end{lem}

For
\begin{align*}
\widetilde{\frak {L}}^{k}u &= {\rm e}^{2y_{1}}\left(-\sum
_{r,s=1}^{2}(a_{r,s}(x)u_{x_{s}})_{x_{r}}+\sum
_{r=1}^{2}b_{r}(x)u_{x_{r}}+c(x)u\right)\\[.2pc]
&= \left(\sum
_{r,s=1}^{2}-(\widetilde{a}_{r,
s}^{k}(y)u_{y_{s}})_{y_{r}}\right) + \left(\sum
_{r=1}^{2}\widetilde{b}_{r}^{k}(y)u_{y_{r}}+\widetilde{c}^{k}
(y)u\right)\\[.2pc]
&=\widetilde{\frak {M}}^{k}u+\widetilde{\frak{N}}^{k}u.
\end{align*}
Here
\begin{equation}
\widetilde{\frak {N}}^{k}u = \sum_{r= 1}^{2}
\widetilde{b}_{r}^{k}(y)u_{y_{r}}+\widetilde{c}^{k} (y)u
\end{equation}
and $y=(y_{1},y_{2})=(\tau _{k},\theta _{k})$ for
some $k$. Moreover the coefficients of $\widetilde{\frak {N}}^{k}$
satisfy
\begin{align*}
\widetilde{b}_{r}^{k} &= O({\rm e}^{\tau _{k}})\; \mathrm{for}\:
r=1,2\\
\intertext{and}
\widetilde{c}^{k} &= O({\rm e}^{2\tau _{k}})
\end{align*}
as $\tau _{k}\rightarrow -\infty .$

Once more
\begin{equation*}
\int _{\widetilde{\Omega }_{i,j}^{k}}|\widetilde{\frak
{M}}^{k}u|^{2}{\rm d}y\leq 2\left(\int _{\widetilde{\Omega
}_{i,j}^{k}}|\widetilde{\frak {L}}^{k}u|^{2}{\rm d}y+\int
_{\widetilde{\Omega }_{i,j}^{k}}|\widetilde{\frak
{N}}^{k}u|^{2}{\rm d}y\right).
\end{equation*}
Using Lemma~\ref{4SLrabsN2u} we can conclude that there exists a
constant $C$ such that the following estimate holds.
\begin{align}
&\frac{\mu _{0}^{2}}{2}\sum _{r,s=1}^{2}\int _{\widetilde{\Omega
}_{i,j}^{k}}\left|\frac{\partial ^{2}u}{\partial y_{r}\partial
y_{s}}\right|^{2}{\rm d}y-C\left(\sum _{r=1}^{2}\int _{\widetilde{\Omega
}_{i,j}^{k}}\left|\frac{\partial u}{\partial y_{r}}\right|^{2}{\rm d}y+\int
_{\widetilde{\Omega }_{i,j}^{k}}|u|^{2}{\rm e}^{4y_{1}}{\rm
d}y\right)\nonumber\\[.2pc]
&\leq 2\int _{\widetilde{\Omega }_{i,j}^{k}}|\widetilde{\frak
{L}}^{k}u|^{2}{\rm d}y+2\sum _{j=1}^{4}\int _{\widetilde{\gamma
}_{j}}\left(\frac{\partial u}{\partial t}\right)_{\widetilde{A}^{k}}\frac{{\rm
d}}{{\rm d}\sigma }\left(\frac{\partial u}{\partial
n}\right)_{\widetilde{A}^{k}}{\rm d}\sigma \nonumber\\[.2pc]
&\quad\ +\sum _{j=1}^{4}\left\{ \left(\frac{\partial u}{\partial
n^{j+1}}\right)_{\widetilde{A}^{k}} \left(\frac{\partial u}{\partial
t^{j+1}}\right)_{\widetilde{A}^{k}}-\left(\frac{\partial u}{\partial
n^{j}}\right)_{\widetilde{A}^{k}}\left(\frac{\partial u}{\partial
t^{j}}\right)_{\widetilde{A}^{k}}\right\} (\widetilde{Q}_{j})\nonumber \\[.2pc]
&\quad\ +\sum _{j=1}^{4}\int _{\widetilde{\gamma
}_{j}}|\widetilde{\kappa }|\left(\left(\frac{\partial u}{\partial
n}\right)_{\widetilde{A}^{k}}^{2}+\left(\frac{\partial u}{\partial
t}\right)_{\widetilde{A}^{k}}^{2}\right){\rm d}\sigma.\label{EabsN012uy}
\end{align}
Rewriting (\ref{EabsN012uy}) in $(\nu _{k},\phi _{k})$
coordinates (\ref{EabsN2unufi}) follows.\hfill $\cd$\vspace{.5pc}

We now need to obtain estimates for the spectral element functions
in the $H^{1}$ norm which we do in the following theorem.

\setcounter{theore}{0}
\begin{theor}[\!]\label{4STuN1} The following estimate holds{\rm :}
\begin{align}
&\sum_{k=1}^{p}\sum _{i=1}^{I_{k}}|u_{i,1}^{k}|^{2}+\sum
_{k=1}^{p}\sum _{j=2}^{M}\sum _{i=1}^{I_{k}}\Vert
u_{i,j}^{k}(\nu _{k},\phi _{k})\Vert
_{1,\widehat{\Omega }_{i,j}^{k}}^{2}+\sum _{l=1}^{L}\Vert
u_{l}^{p+1}(\xi ,\eta )\Vert
_{1,S}^{2}\nonumber\\[.2pc]
&\leq C_{M}\left\{ \sum _{k=1}^{p}\sum _{j=2}^{M}\sum
_{i=1}^{I_{k}}\Vert \frak {L}_{i,j}^{k}u_{i,j}^{k}(\nu
_{k},\phi _{k})\Vert _{0,\widehat{\Omega
}_{i,j}^{k}}^{2}\begin{array}{c}\\ \\ \\[-.15pc] \end{array}\right.\nonumber \\[.2pc]
&\quad\ +\sum _{k=1}^{p}\sum _{\gamma _{s}\subseteq \Omega
^{k}}(\Vert [u]\Vert _{0,\widehat{\gamma
}_{s}}^{2}+\Vert [u_{\nu _{k}}]\Vert
_{0,\widehat{\gamma }_{s}}^{2}+\Vert [u_{\phi
_{k}}]\Vert _{0,\widehat{\gamma }_{s}}^{2})\nonumber\\[.2pc]
&\quad\ +\sum _{l\in \mathcal{D}}\sum _{k=l-1}^{l}\sum _{\gamma
_{s}\subseteq \partial \Omega ^{k}\bigcap \Gamma _{l},\mu
(\widehat{\gamma }_{s})<\infty }(\Vert u\Vert
_{0,\widehat{\gamma }_{s}}^{2}+\Vert u_{\nu _{k}}\Vert
_{0,\widehat{\gamma }_{s}}^{2})\nonumber
\end{align}
\begin{align}
&\quad\ +\sum _{k=1}^{p}\sum _{\gamma _{s}\subseteq B_{\rho
}^{k}}(\Vert [u]\Vert _{0,\widehat{\gamma
}_{s}}^{2}+\Vert [u_{\nu _{k}}]\Vert
_{0,\widehat{\gamma }_{s}}^{2}+\Vert [u_{\phi
_{k}}]\Vert _{0,\widehat{\gamma }_{s}}^{2})\nonumber\\[.2pc]
&\quad\ +\sum _{l\in \mathcal{N}}\sum _{k=l-1}^{l}\sum _{\gamma
_{s}\subseteq \partial \Omega ^{k}\bigcap \Gamma _{l}}\left\Vert
\left(\frac{\partial u}{\partial
n}\right)_{\widetilde{A}^{k}}\right\Vert _{0,\widehat{\gamma
}_{s}}^{2}\nonumber \\[.2pc]
&\quad\ +\sum _{l=1}^{L}\int _{S}\int |\frak
{L}_{l}^{p+1}u_{l}^{p+1}(\xi ,\eta )|^{2}{\rm d}\xi {\rm d}\eta\nonumber\\[.2pc]
&\quad\ +\sum _{\gamma _{s}\subseteq \Omega ^{p+1}}(\Vert
[u]\Vert _{0,\gamma _{s}}^{2}+\Vert
[u_{x_{1}}]\Vert _{0,\gamma _{s}}^{2}+\Vert
[u_{x_{2}}]\Vert _{0,\gamma _{s}}^{2})\nonumber\\[.2pc]
&\quad\ +\sum _{l\in \mathcal{D}}\sum _{\gamma _{s}\subseteq \partial
\Omega ^{p+1}\bigcap \Gamma _{l}}\left(\Vert u\Vert
_{0,\gamma _{s}}^{2}+\left\Vert \frac{\partial u}{\partial T}\right\Vert
_{0,\gamma _{s}}^{2}\right)\nonumber\\[.2pc]
&\quad\ + \left.\sum _{l\in \mathcal{N}}\sum _{\gamma
_{s}\subseteq \partial \Omega ^{p+1}\bigcap \Gamma _{l}}\left\Vert
\left(\frac{\partial u}{\partial N}\right)_{A}\right\Vert _{0,\gamma
_{s}}^{2}\right\}.\label{4SETuN1}
\end{align}
Here $C_{M}=CM^{4}$ if there exists a vertex $A_{j}$ such that Neumann
boundary conditions are imposed on the adjoining sides $\Gamma _{j}$ and
$\Gamma _{j+1}$ and $C_{M}=C$ otherwise. $C$ denotes a constant and $\mu
(\hat{{\gamma }_{s}})$ the length of $\widehat{\gamma}_{s}$.
\end{theor}

To prove the estimate (\ref{4SETuN1}) we shall use (\ref{Eapriori}).
To do so we have to define a corrected version of the spectral element
functions so that it is conforming.

Let $\{ \{ u_{i,j}^{k}(\nu _{k},\phi _{k})\}
_{i,j\leq M,k},\{ u_{i,j}^{k}(\xi ,\eta )\}
_{i,j>M,k}\} _{k}$ be a set of spectral element functions $\in \pi
^{M,W}.$ Here $\pi ^{M,W}$ is the set of spectral element functions such
that $u_{i,1}^{k}=g_{k}$, a constant for all $i$, and $u_{i,j}^{k}$ is a
polynomial of degree $W$ in each variable for $j\geq 2$. Then there is a
set of spectral element functions
\begin{equation*}
\{ \lambda _{i,j}^{k}(\nu _{k},\phi _{k})\}
_{i,j\leq M,k},\{ \lambda _{i,j}^{k}(\xi ,\eta )\}
_{i,j>M,k}\in \pi ^{M,W}
\end{equation*}
such that the function $\varphi (x_{1},x_{2})$ defined as
\begin{align*}
&\varphi (x_{1},x_{2})\\[.2pc]
&=\left\{ \begin{array}{c}
(u_{i,j}^{k}+\lambda _{i,j}^{k})(\nu
_{k}(x_{1},x_{2}),\phi _{k}(x_{1},x_{2}))\ \
\hbox{if}\ \ (x_{1},x_{2})\in \Omega _{i,j}^{k}\ \ \hbox{for}\ \ j\leq M\\[.3pc]
(u_{i,j}^{k}+\lambda _{i,j}^{k})(\xi
(x_{1},x_{2}),\eta (x_{1},x_{2}))\: \ \ \hbox{if}\ \
(x_{1},x_{2})\in \Omega _{i,j}^{k}\ \ \hbox{for}\ \ j > M
\end{array} \right.
\end{align*}
is a differentiable function of its arguments and $\varphi \in
H_{0}^{1}(\Omega ).$ This can be shown as in Lemma~4.57 of
\cite{schwab}.

Moreover the estimate
\begin{align}
&\sum _{k=1}^{p}\sum _{i=1}^{I_{k}}|\lambda
_{i,1}^{k}|^{2}+\sum _{k=1}^{p}\sum _{j=2}^{M}\sum
_{i=1}^{I_{k}}\Vert \lambda _{i,j}^{k}(\nu _{k},\phi
_{k})\Vert _{1,\widehat{\Omega }_{i,j}^{k}}^{2}\nonumber\\[.2pc]
&\quad\ +\sum_{k=1}^{p}\sum _{j=M+1}^{J_{k}}\sum _{i=1}^{I_{k,j}}\Vert \lambda
_{i,j}^{k}(\xi ,\eta )\Vert _{1,S}^{2}\nonumber
\end{align}
\begin{align}
&\leq C\left\{ \left(\sum _{l\in \mathcal{D}}\sum _{k=l-1}^{l}\sum
_{\gamma _{s}\subseteq \Gamma _{l}\cap \partial \Omega ^{k},\mu
(\widehat{\gamma }_{s})<\infty }(\Vert u\Vert
_{0,\widehat{\gamma }_{s}}^{2}+\Vert u_{\nu _{k}}\Vert
_{0,\widehat{\gamma }_{s}}^{2})\right)\right.\nonumber \\[.2pc]
&\quad\ +\sum _{k=1}^{p}\sum _{\gamma _{s}\subseteq \Omega ^{k},\mu
(\widehat{\gamma }_{s})<\infty }(\Vert
[u]\Vert _{0,\widehat{\gamma }_{s}}^{2}+\Vert
[u_{\nu _{k}}]\Vert _{0,\widehat{\gamma
}_{s}}^{2}+\Vert [u_{\phi _{k}}]\Vert
_{0,\widehat{\gamma }_{s}}^{2})\nonumber \\[.2pc]
&\quad\ +\sum _{k=1}^{p}\sum _{\gamma _{s}\subseteq B_{\rho
}^{k}}(\Vert [u]\Vert _{0,\widehat{\gamma
}_{s}}^{2}+\Vert [u_{\nu _{k}}]\Vert
_{0,\widehat{\gamma }_{s}}^{2}+\Vert [u_{\phi
_{k}}]\Vert _{0,\widehat{\gamma }_{s}}^{2})\nonumber\\[.2pc]
&\quad\ + \sum _{\gamma _{s}\subseteq \Omega
^{p+1}}\left(\Vert [u]\Vert _{0,\gamma
_{s}}^{2}+\left\Vert \left[\frac{\partial u}{\partial
T}\right]\right\Vert _{0,\gamma _{s}}^{2}\right)\nonumber \\[.2pc]
&\quad\ \left. +\sum _{l\in
\mathcal{D}}\sum _{\gamma _{s}\subseteq \partial \Omega ^{p+1}\cap
\Gamma _{l}}\left(\Vert u\Vert _{0,\gamma
_{s}}^{2}+\left\Vert \frac{\partial u}{\partial T}\right\Vert _{0,\gamma
_{s}}^{2}\right)\right\}\label{EN1lamall}
\end{align}
holds.

We now explain the notation we have used in (\ref{EN1lamall}). Let ${\rm
d}\widehat{\sigma }$ denote an element of arc length in $(\nu
_{k},\phi _{k})$ coordinates. Then
\begin{equation*}
\Vert w\Vert _{0,\widehat{\gamma }_{s}}^{2}=\int
_{\widehat{\gamma }_{s}}|w(\nu _{k},\phi
_{k})|^{2}{\rm d}\widehat{\sigma }.
\end{equation*}
Moreover if $\gamma _{s}$ is given by $\gamma _{s}=\partial \Omega
_{m}^{p+1}\bigcap \partial \Omega _{n}^{p+1}$ then
\begin{equation*}
\left\Vert \left[\frac{\partial u}{\partial T}\right]\right\Vert
_{0,\gamma _{s}}^{2}=\int _{\gamma _{s}}\left(\frac{\partial
u_{m}^{p+1}}{\partial T}-\frac{\partial u_{n}^{p+1}}{\partial
T}\right)^{2}{\rm d}s.
\end{equation*}
Here ${\partial}/{\partial T}$ denotes the tangential derivative in
$(x_{1},x_{2})$ variables, i.e.
\begin{equation*}
\frac{\partial u}{\partial T}=T^{t}\nabla _{x}u.
\end{equation*}
The other terms in the right-hand side of (\ref{EN1lamall}) are
similarly defined.

Now consider the bilinear form
\begin{align*}
B(\varphi ,v) &= \int _{\Omega }\left(\sum
_{r,s=1}^{2}a_{r,s}(x)\varphi _{x_{s}}v_{x_{r}}+\sum
_{r=1}^{2}b_{r}(x)\varphi _{x_{r}}v+c\varphi v\right){\rm d}x\\[.2pc]
&= \sum _{k=1}^{p}\sum _{j=1}^{M}\sum _{i=1}^{I_{k}}B(\varphi
,v)_{\Omega _{i,j}^{k}}+\sum _{l=1}^{L}B(\varphi
,v)_{\Omega _{l}^{p+1}}.
\end{align*}
Here
\begin{equation*}
B(\varphi ,v)_{\Delta }=\int _{\Delta }\left(\sum
_{r,s=1}^{2}a_{r,s}(x)\varphi _{x_{s}}v_{x_{r}}+\sum
_{r=1}^{2}b_{r}(x)\varphi _{x_{r}}v+c\varphi
v\right){\rm d}x,
\end{equation*}
where $\Delta$ is a domain contained in $\Omega $ and $v\in
H_{0}^{1}(\Omega )$.

Now
\begin{align*}
B(\varphi ,v)_{\Omega _{l}^{p+1}} &= \int _{\Omega
_{l}^{p+1}}\left(\sum _{r,s=1}^{2}a_{r,s}(x)\varphi
_{x_{s}}v_{x_{r}}+\sum _{r=1}^{2}b_{r}(x)\varphi
_{x_{r}}v+c\varphi v\right){\rm d}x\\[.2pc]
&= \int _{\Omega _{l}^{p+1}}\mathfrak{L}\varphi v{\rm d} x + \int
_{\partial \Omega _{l}^{p+1}}\left(\frac{\partial \varphi }{\partial
N}\right)_{A}v{\rm d}s.
\end{align*}
Similarly if $1\leq j\leq M$ we have
\begin{equation*}
B(\varphi ,v)_{\Omega _{i,j}^{k}}=\int _{\widetilde{\Omega
}_{i,j}^{k}}\widetilde{\frak {L}}^{k}\varphi v{\rm d}\tau _{k}{\rm d}\theta
_{k}+\int _{\partial \widetilde{\Omega }_{i,j}^{k}}\left(\frac{\partial
\varphi }{\partial n}\right)_{\widetilde{A}^{k}}v{\rm d}\sigma.
\end{equation*}
Moreover if $j=1$,
\begin{equation*}
B(\varphi ,v)_{\Omega _{i,1}^{k}}=\int _{\widetilde{\Omega
}_{i,1}^{k}}c\varphi v{\rm e}^{2\tau _{k}}{\rm d}\tau _{k}{\rm d}\theta
_{k}+\int _{\partial \widetilde{\Omega }_{i,1}^{k}}\left(\frac{\partial
\varphi }{\partial n}\right)_{\widetilde{A}^{k}}v{\rm d}\sigma
\end{equation*}
since $\varphi $ is a constant on $\widetilde{\Omega }_{i,1}^{k}$.

Finally if $j=M+1$ we obtain
\begin{align*}
B(\varphi ,v)_{\Omega _{i,M+1}^{k}} &= \int _{\Omega
_{i,M+1}^{k}}\mathfrak{L}\varphi v{\rm d}x + \int _{\widetilde{B}_{\rho
}^{k}}\left(\frac{\partial \varphi }{\partial
n}\right)_{\widetilde{A}^{k}}v{\rm d}\sigma\\[.2pc]
&\quad\ + \int _{\partial \Omega _{i,M+1}^{k}\setminus B_{\rho
}^{k}}\left(\frac{\partial \varphi }{\partial N}\right)_{A}v{\rm d}s.
\end{align*}
For by (\ref{EuNAunAtil})
\begin{equation*}
\rho \left(\frac{\partial \varphi }{\partial
N}\right)_{A}(P) = \left(\frac{\partial \varphi }{\partial
n}\right)_{\widetilde{A}^{k}}(\widetilde{P})
\end{equation*}
and ${\rm d}s=\rho {\rm d}\sigma$. Here $P$ is any point on the circular arc
$B_{\rho }^{k}$ and $\widetilde{P}$ is its image in $(\tau
_{k},\theta _{k})$ coordinates. Now
\begin{align}
B(\varphi ,v) &= \sum _{k=1}^{p}\sum _{j=1}^{M}\sum
_{i=1}^{I_{k}}B(\varphi ,v)_{\Omega _{i,j}^{k}}+\sum
_{l=1}^{L}B(\varphi ,v)_{\Omega _{l}^{p+1}}\nonumber \\[.2pc]
&= \sum _{k=1}^{p}\sum _{j=1}^{M}\sum
_{i=1}^{I_{k}}B(u_{i,j}^{k},v)_{\Omega _{i,j}^{k}}+\sum
_{l=1}^{L}B(u_{l}^{p+1},v)_{\Omega
_{l}^{p+1}}\nonumber\\[.2pc]
&\quad\ +\left(\sum _{k=1}^{p}\sum _{j=1}^{M}\sum
_{i=1}^{I_{k}}B(\lambda _{i,j}^{k},v)_{\Omega
_{i,j}^{k}}+\sum _{l=1}^{L}B(\lambda _{l}^{p+1},v)_{\Omega
_{l}^{p+1}}\right)\nonumber\\[.2pc]
&= \sum _{k=1}^{p}\sum _{j=1}^{M}\sum _{i=1}^{I_{k}}\int
_{\widetilde{\Omega }_{i,j}^{k}}\widetilde{\frak
{L}}^{k}u_{i,j}^{k}v{\rm d}\tau _{k}{\rm d}\theta _{k}+\sum
_{l=1}^{L}\int _{\Omega _{l}^{p+1}}\mathfrak{L}u_{l}^{p+1}v{\rm
d}x_{1}{\rm d}x_{2}\nonumber \\[.2pc]
&\quad\ +\sum _{k=1}^{p}\sum _{\gamma _{s}\subseteq \Omega ^{k},\mu
(\widetilde{\gamma }_{s})<\infty }\int _{\widetilde{\gamma
}_{s}}\left[\left(\frac{\partial u}{\partial
n}\right)_{\widetilde{A}^{k}}\right]v{\rm d}\sigma \nonumber
\end{align}
\begin{align}
&\quad\ +\sum _{\gamma _{s}\subseteq \Omega ^{p+1}}\int _{\gamma
_{s}}\left[\left(\frac{\partial u}{\partial N}\right)_{A}\right]v{\rm
d}s+\sum _{k=1}^{p}\sum _{\gamma _{s}\subseteq B_{\rho }^{k}}\int
_{\widetilde{\gamma }_{s}}\left[\left(\frac{\partial u}{\partial
n}\right)_{\widetilde{A}^{k}}\right]v{\rm d}\sigma \nonumber \\[.2pc]
&\quad\ +\sum _{l\in \mathcal{N}}\sum _{k=l-1}^{l}\sum _{\gamma
_{s}\subseteq \Gamma _{l}\cap \partial \Omega ^{k},\mu
(\widetilde{\gamma }_{s})<\infty }\int _{\widetilde{\gamma
}_{s}}\left(\frac{\partial u}{\partial
n}\right)_{\widetilde{A}^{k}}v{\rm d}\sigma\nonumber\\[.2pc]
&\quad\ +\sum _{l\in \mathcal{N}}\sum _{r=1}^{L}\sum _{\gamma
_{s}\subseteq \partial \Omega _{r}^{p+1}\cap \Gamma _{l}}\int _{\gamma
_{s}}\left(\frac{\partial u}{\partial N}\right)_{A}v{\rm d}s\nonumber \\[.2pc]
&\quad\ + \left(\sum _{k=1}^{p}\sum _{i=1}^{I_{k}}B(\lambda
_{i,1}^{k},v)_{\Omega _{i,1}^{k}}+\sum _{k=1}^{p}\sum _{j=2}^{M}\sum
_{i=1}^{I_{k}}B(\lambda _{i,j}^{k},v)_{\Omega _{i,j}^{k}}
\right.\nonumber\\[.2pc]
&\quad\ \left. +\sum _{l=1}^{L}B(\lambda _{l}^{p+1},v)_{\Omega
_{l}^{p+1}}\right).\label{Ebilin}
\end{align}
Now
\begin{equation*}
\int _{\widetilde{\Omega }_{i,1}^{k}}\widetilde{\frak {L}}^{k}\lambda
_{i,1}^{k}v{\rm d}\tau _{k}{\rm d}\theta _{k}=\int _{\widetilde{\Omega
}_{i,1}^{k}}c\lambda _{i,1}^{k}v{\rm e}^{2\tau _{k}}{\rm d}\tau _{k}{\rm
d}\theta _{k}.
\end{equation*}
Here
\begin{equation*}
\lambda _{i,1}^{k} = \left\{ \begin{array}{cl}
-u_{i,1}^{k}, & \textrm{ if }\Gamma _{k}\textrm{ or }\Gamma
_{k+1}\subseteq \Gamma ^{[0]}\\[.3pc]
 0, & \textrm{otherwise}\end{array}
\right\}.
\end{equation*}
Now $c_{k}=c(A_{k}),$ a constant, and $c(x_{1},x_{2})$ is an analytic
function of $x_{1}$ and $x_{2}$. Hence
\begin{align*}
\left|\int _{\widetilde{\Omega }_{i,1}^{k}}\widetilde{\frak
{L}}^{k}\lambda _{i,1}^{k}v{\rm d}\tau _{k}{\rm d}\theta _{k}\right| &\leq
2c_{k}\left(\int |\lambda _{i,1}^{k}|^{2}{\rm e}^{2\tau _{k}}{\rm d}\tau
_{k}{\rm d}\theta _{k}\right)^{1/2}\\[.2pc]
&\quad\ \times \left(\int v^{2}{\rm e}^{2\tau _{k}}{\rm d}\tau
_{k}{\rm d}\theta _{k}\right)^{1/2}
\end{align*}
for $M$ large enough. And so we obtain
\begin{equation*}
\left|\int _{\widetilde{\Omega }_{i,1}^{k}}\widetilde{\frak
{L}}^{k}\lambda _{i,1}^{k}v{\rm d}\tau _{k}{\rm d}\theta _{k}\right|\leq
\varepsilon |\lambda _{i,1}^{k}|\Vert v(x_{1},x_{2})\Vert _{0,\Omega
_{i,1}^{k}},
\end{equation*}
where $\varepsilon $ is exponentially small in $M.$ Now, let $2\leq
j\leq M.$ Then
\begin{equation*}
\left|\int _{\widetilde{\Omega }_{i,j}^{k}}\widetilde{\frak
{L}}^{k}u_{i,j}^{k}v{\rm d}\tau _{k}{\rm d}\theta _{k}\right| \leq
\Vert \widetilde{\frak {L}}^{k}u_{i,j}^{k}(\tau _{k},\theta
_{k})\Vert _{0,\widetilde{\Omega }_{i,j}^{k}}\Vert
v(\tau _{k},\theta _{k})\Vert _{0,\widetilde{\Omega}_{i,j}^{k}}.
\end{equation*}
Finally
\begin{equation*}
\left|\int _{\Omega_{l}^{p+1}}(\mathfrak{L}u_{l}^{p+1})v{\rm
d}x\right|\leq \Vert \mathfrak{L}u_{l}^{p+1}(x_{1},x_{2})\Vert
_{0,\Omega _{l}^{p+1}}\Vert v(x_{1},x_{2})\Vert _{0,\Omega _{l}^{p+1}}.
\end{equation*}
Now
\begin{equation*}
\sum _{k=1}^{p}\sum _{j=2}^{M}\sum _{i=1}^{I_{k}}\Vert v(\nu _{k},\phi
_{k})\Vert _{0,\widehat{\Omega }_{i,j}^{k}}^{2}\leq K_{M}\Vert
v(x_{1},x_{2})\Vert _{1,\Omega }^{2}.
\end{equation*}
Here $K_{M}=KM^{2}$ if there is a vertex $A_{j}$ such that Neumann
boundary conditions are imposed on the adjoining sides $\Gamma _{j}$ and
$\Gamma _{j+1}$ and $K_{M}=K$, otherwise. $K$ denotes a constant. Hence
\begin{equation}
\sum _{k=1}^{p}\sum _{j=2}^{M}\sum _{i=1}^{I_{k}}\Vert v(\nu _{k},\phi
_{k})\Vert _{1,\widehat{\Omega }_{i,j}^{k}}^{2}\leq K_{M}\Vert
v(x_{1},x_{2})\Vert _{1,\Omega }^{2}.\label{EN1votilo}
\end{equation}
Now using the trace theorem for Sobolev spaces we obtain
\begin{equation*}
\sum _{k=1}^{p}\sum _{j=2}^{M}\sum _{i=1}^{I_{k}}\Vert v\Vert
_{0,\partial \widehat{\Omega }_{i,j}^{k}}^{2}\leq K_{M}\Vert
v(x_{1},x_{2})\Vert _{1,\Omega }^{2}.
\end{equation*}
And so we can conclude that
\begin{equation}
\sum _{k=1}^{p}\sum _{j=2}^{M}\sum _{i=1}^{I_{k}}\int _{\partial
\widetilde{\Omega }_{i,j}^{k}}v^{2}{\rm d}\sigma \leq K_{M}\Vert
v(x_{1},x_{2})\Vert _{1,\Omega
}^{2}.\label{ENL2vN1v}
\end{equation}
Using the Cauchy--Schwartz inequality in (\ref{Ebilin}) and using
(\ref{EN1votilo}) and (\ref{ENL2vN1v}) we can conclude that
\begin{align*}
|B(\varphi ,v)|^{2} &\leq K\left\{ \sum
_{k=1}^{p}\sum _{j=2}^{M}\sum _{i=1}^{I_{k}}\Vert \widetilde{\frak
{L}}^{k}u_{i,j}^{k}(\tau _{k},\theta _{k})\Vert
_{0,\widetilde{\Omega }_{i,j}^{k}}^{2}+\sum _{k=1}^{p}\sum
_{i=1}^{I_{k}}\varepsilon |u_{i,1}^{k}|^{2}\right.\\[.2pc]
&\quad\ +\sum _{k=1}^{p} \left( \begin{array}{@{}c@{}}\\ \\ \\[-.3pc] \end{array} \right. \sum _{\gamma _{s}\subseteq \Omega
^{k}}\int _{\widetilde{\gamma }_{s}}\left[\left(\!\frac{\partial
u}{\partial n}\!\right)_{\widetilde{A}^{k}}\right]^{2}{\rm d}\sigma \!+\!\sum
_{\gamma _{s}\subseteq B_{\rho }^{k}}\int _{\widetilde{\gamma
}_{s}}\left[\left(\!\frac{\partial u}{\partial
n}\!\right)_{\widetilde{A}^{k}}\right]^{2}{\rm d}\sigma \left. \begin{array}{@{}c@{\!\!}}\\ \\ \\[-.3pc] \end{array} \right)\\[.2pc]
&\quad\ +\sum _{l\in \mathcal{N}}\sum _{k=l-1}^{l}\sum _{\gamma _{s}\subseteq
\Gamma _{l}\bigcap \partial \Omega ^{k}}\int _{\widetilde{\gamma
}_{s}}\left(\frac{\partial u}{\partial
n}\right)_{\widetilde{A}^{k}}^{2}{\rm d}\sigma\\[.2pc]
&\quad\ +\sum _{l=1}^{L}\int _{\Omega _{l}^{p+1}}\int
|\mathfrak{L}u_{l}^{p+1}(x_{1},
x_{2})|^{2}{\rm d}x_{1}{\rm d}x_{2}\\[.2pc]
&\quad\ +\sum _{\gamma _{s}\subseteq \Omega
^{p+1}}\int _{\gamma _{s}}\left[\left(\frac{\partial u}{\partial
N}\right)_{A}\right]^{2}{\rm d}s\\[.2pc]
&\quad\ +\sum _{l\in \mathcal{N}}\sum _{\gamma _{s}\subseteq \Gamma
_{l}\bigcap \partial \Omega ^{p+1}}\int _{\gamma
_{s}}\left(\frac{\partial u}{\partial
N}\right)_{A}^{2}{\rm d}s\\[.2pc]
&\quad\ +\varepsilon \sum _{k=1}^{p}\sum
_{i=1}^{I_{k}}|\lambda _{i,1}^{k}|^{2}+\sum _{k=1}^{p}\sum
_{j=2}^{M}\sum _{i=1}^{I_{k}}\Vert \lambda _{i,j}^{k}(\tau
_{k},\theta _{k})\Vert _{1,\widetilde{\Omega
}_{i,j}^{k}}^{2}
\end{align*}
\begin{align*}
&\quad\ \left. + \sum _{l=1}^{L}\Vert \lambda
_{l}^{p+1}(x,y)\Vert _{1,\Omega _{l}^{p+1}}^{2} \right\}
\cdot \left\{ \sum _{k=1}^{p}\sum _{i=1}^{I_{k}}\Vert
v(x_{1},x_{2})\Vert _{0,\Omega _{i,1}^{k}}^{2} \right.\\[.2pc]
&\quad\ +\sum _{k=1}^{p}\sum _{j=2}^{M}\sum _{i=1}^{I_{k}}\Vert
v(\tau _{k},\theta _{k})\Vert _{1,\widetilde{\Omega
}_{i,j}^{k}}^{2}\\[.2pc]
&\quad\ + \sum _{l=1}^{L}\Vert
v_{l}^{p+1}(x,y)\Vert _{1,\Omega _{l}^{p+1}}^{2}+\sum
_{k=1}^{p}\sum _{j=2}^{M}\sum _{i=1}^{I_{k}}\int _{\partial
\widetilde{\Omega }_{i,j}^{k}}v^{2}{\rm d}\sigma\\[.2pc]
&\quad\ + \sum _{\gamma _{s}\subseteq \Omega ^{p+1}}\int _{\gamma
_{s}}v^{2}{\rm d}s+\sum _{l\in \mathcal{N}}\sum _{\gamma _{s}\subseteq
\partial \Omega ^{p+1}\bigcap \Gamma _{l}}\int _{\gamma
_{s}}v^{2}{\rm d}s \left. \begin{array}{@{}c@{}}\\ \\
\\[-.3pc] \end{array} \right\}.
\end{align*}
Now $v\in H_{0}^{1}(\Omega )$ and $\mathfrak{L}$ satisfies
the inf--sup conditions (3.4). Hence using (\ref{Eapriori}),
(\ref{EN1lamall}) and (\ref{ENL2vN1v}) we obtain
\begin{align*}
\Vert \varphi \Vert _{1,\Omega }^{2} &\leq K_{M} \left\{ \sum
_{k=1}^{p}\sum _{j=2}^{M}\sum _{i=1}^{I_{k}}\Vert \widetilde{\frak
{L}}^{k}u_{i,j}^{k}(\tau _{k},\theta _{k})\Vert
_{0,\widetilde{\Omega }_{i,j}^{k}}^{2}\right.\\[.2pc]
&\quad\ +\sum _{k=1}^{p} \left( \begin{array}{@{}c@{}}\\ \\
\\[-.3pc] \end{array} \right. \sum _{\gamma _{s}\subseteq \Omega
^{k}}(\Vert [u]\Vert _{0,\widehat{\gamma
}_{s}}^{2}+\Vert [u_{\nu _{k}}]\Vert
_{0,\widehat{\gamma }_{s}}^{2}+\Vert [u_{\phi
_{k}}]\Vert _{0,\widehat{\gamma }_{s}}^{2})\\[.2pc]
&\quad\ +\sum _{l\in \mathcal{N}}\sum _{k=l-1}^{l}\sum _{\gamma
_{s}\subseteq \partial \Omega ^{k}\bigcap \Gamma _{l}}\int
_{\widetilde{\gamma }_{s}}\left(\frac{\partial u}{\partial
n}\right)_{\widetilde{A}^{k}}^{2}{\rm d}\sigma \\[.2pc]
&\quad\ + \sum _{l\in \mathcal{D}}\sum _{k=l-1}^{l}\sum _{\gamma
_{s}\subseteq \Gamma _{l}\bigcap \partial \Omega ^{k},\mu
(\widehat{\gamma }_{s})<\infty }(\Vert u\Vert
_{0,\widehat{\gamma }_{s}}^{2}+\Vert u_{\nu _{k}}\Vert
_{0,\widehat{\gamma }_{s}}^{2}) \left. \begin{array}{@{}c@{}}\\ \\
\\[-.3pc] \end{array} \right)\\[.2pc]
&\quad\ +\sum _{k=1}^{p}\sum _{\gamma _{s}\subseteq B_{\rho
}^{k}}(\Vert [u]\Vert _{0,\widehat{\gamma
}_{s}}^{2}+\Vert [u_{\nu _{k}}]\Vert
_{0,\widehat{\gamma }_{s}}^{2}+\Vert [u_{\phi
_{k}}]\Vert _{0,\widehat{\gamma }_{s}}^{2})\\[.2pc]
&\quad\ +\sum _{l=1}^{L}\int _{\Omega _{l}^{p+1}}\int |\frak
{L}u_{l}^{p+1}(x_{1},x_{2})|^{2}{\rm d}x_{1}{\rm
d}x_{2}\\[.2pc]
&\quad\ +\sum _{\gamma _{s}\subseteq \Omega ^{p+1}}(\Vert
[u]\Vert _{0,\gamma _{s}}^{2}+\Vert
[u_{x_{1}}]\Vert _{0,\gamma _{s}}^{2}+\Vert
[u_{x_{2}}]\Vert _{0,\gamma _{s}}^{2})\\[.2pc]
&\quad\ +\sum _{l\in \mathcal{D}}\sum _{\gamma _{s}\subseteq \partial
\Omega ^{p+1}\bigcap \Gamma _{l}}\left(\Vert u\Vert
_{0,\gamma _{s}}^{2}+\left\Vert \frac{\partial u}{\partial T}\right\Vert
_{0,\gamma _{s}}^{2}\right)\\[.2pc]
&\quad\ + \sum _{l\in \mathcal{N}}\sum _{\gamma _{s}\subseteq
\partial \Omega ^{p+1}\bigcap \Gamma _{l}}\int \left(\frac{\partial
u}{\partial N}\right)_{A}^{2}{\rm d}s
\left. +\varepsilon \sum _{k=1}^{p}\sum
_{i=1}^{I_{k}}(|u_{i,1}^{k}|^{2}+|\lambda
_{i,1}^{k}|^{2})\right\}.
\end{align*}
Here $\varepsilon $ is exponentially small in $M$.

Using (\ref{EN1lamall}) and (\ref{EN1votilo}) once more we obtain
the result.\hfill $\cd$\vspace{.5pc}

\pagebreak

We now define differential operators $(\mathfrak{L}_{i,j}^{k})^{a}$
which are second order differential operators with polynomial
coefficients in $\nu _{k}$ and $\phi_{k}$ of degree $W$ such that these
coefficients are exponentially close approximation to the coefficients
of $(\mathfrak{L}_{i,j}^{k})$ as has been described in the beginning of
this section. In the same way we define the differential operator
$({\partial u}/{\partial n})_{\widetilde{A}^{k}}^{a}$ to be a first
order differential operator with polynomial coefficients in $\nu _{k}$
and $\phi _{k}$ such that these coefficients are exponentially close
approximations to the coefficients of $({\partial u}/{\partial
n})_{\widetilde{A}^{k}}$. The other approximations are similarly
defined.

From the above, it is easy to conclude that
\begin{align}
&\sum _{k=1}^{p}\sum _{i=1}^{I_{k}}\left(|u_{i,1}^{k}|^{2}+\sum
_{j=2}^{M}\Vert u_{i,j}^{k}(\nu _{k},\phi _{k})\Vert _{1,\widehat{\Omega
}_{i,j}^{k}}^{2}\right)\nonumber\\[.2pc]
&\quad\ +\sum _{l=1}^{L}\Vert u_{l}^{p+1}(\xi ,\eta )\Vert
_{1,S}^{2}\leq C_{M}(\mathcal{I}),\label{EN1allurhsi}
\end{align}
where
\begin{align*}
\mathcal{I} &= \left\{ \sum _{k=1}^{p}\sum _{j=2}^{M}\sum _{i=1}^{I_{k}}\Vert
(\mathfrak{L}_{i,j}^{k})^{a}u_{i,j}^{k}(\nu _{k},\phi _{k})\Vert
_{0,\widehat{\Omega }_{i,j}^{k}}^{2} \begin{array}{c}\\ \\ \\[-.3pc] \end{array} \right.\\[.2pc]
&\quad\ +\sum _{k=1}^{p}\sum _{\gamma _{s}\subseteq \Omega ^{k}}(\Vert
[u]\Vert _{0,\widehat{\gamma }_{s}}^{2}+\Vert [u_{\nu _{k}}]\Vert
_{0,\widehat{\gamma }_{s}}^{2}+\Vert [u_{\phi _{k}}]\Vert
_{0,\widehat{\gamma }_{s}}^{2})\\[.2pc]
&\quad\ +\sum _{l\in \mathcal{D}}\sum _{k=l-1}^{l}\sum _{\gamma
_{s}\subseteq \partial \Omega ^{k}\bigcap \Gamma _{l},\mu
(\widehat{\gamma }_{s})<\infty }(\Vert u\Vert
_{0,\widehat{\gamma }_{s}}^{2}+\Vert u_{\nu _{k}}\Vert
_{0,\widehat{\gamma }_{s}}^{2})\\[.2pc]
&\quad\ +\sum _{l\in \mathcal{N}}\sum _{k=l-1}^{l}\sum _{\gamma
_{s}\subseteq \partial \Omega ^{k}\bigcap \Gamma _{l},\mu
(\widehat{\gamma }_{s})<\infty }\left\Vert \left(\frac{\partial u}{\partial
n}\right)_{\widetilde{A}^{k}}^{a}\right\Vert _{0,\widehat{\gamma
}_{s}}^{2}\\[.2pc]
&\quad\ +\sum _{k=1}^{p}\sum _{\gamma _{s}\subseteq B_{\rho }^{k}}(\Vert
[u]\Vert _{0,\widehat{\gamma }_{s}}^{2}+\Vert [(u)_{\nu _{k}}^{a}]\Vert
_{0,\widehat{\gamma }_{s}}^{2}+\Vert [(u)_{\phi _{k}}^{a}]\Vert
_{0,\widehat{\gamma }_{s}}^{2})\\[.2pc]
&\quad\ +\sum _{l=1}^{L}\Vert
(\mathfrak{L}_{l}^{p+1})^{a}u_{l}^{p+1}(\xi ,\eta )\Vert
_{0,S}^{2}\\[.2pc]
&\quad\ +\sum _{\gamma _{s}\subseteq \Omega
^{p+1}}(\Vert [u]\Vert _{0,\gamma
_{s}}^{2}+\Vert [u_{x_{1}}]^{a}\Vert _{0,\gamma
_{s}}^{2}+\Vert [u_{x_{2}}]^{a}\Vert _{0,\gamma
_{s}}^{2})\\[.2pc]
&\quad\ +\sum _{l\in \mathcal{D}}\sum _{\gamma _{s}\subseteq \partial
\Omega ^{p+1}\bigcap \Gamma _{l}}\left(\Vert u\Vert
_{0,\gamma _{s}}^{2}+\left\Vert \left(\frac{\partial u}{\partial
T}\right)^{a}\right\Vert _{0,\gamma _{s}}^{2}\right)\\[.2pc]
&\quad\ + \sum _{l\in
\mathcal{N}}\sum _{\gamma _{s}\subseteq \partial \Omega ^{p+1}\bigcap
\Gamma _{l}}\left\Vert \left(\frac{\partial u}{\partial
N}\right)_{A}^{a}\right\Vert _{0,\gamma _{s}}^{2} \left. \begin{array}{@{}c@{}}\\ \\ \\[-.3pc] \end{array} \right\}.
\end{align*}
Here $C_{M}$ is as defined in Theorem \ref{4STuN1}.

\subsection{\it \label{sec3.3:Gen-stab-est}The estimates}

We now define the quadratic form
\begin{align}
&\mathcal{V}^{M,W}(\{ u_{i,j}^{k}(\nu _{k},\phi _{k})\} _{i,j,k},\{
u_{l}^{p+1}(\xi ,\eta )\} _{l})\nonumber\\[.2pc]
&= \left\{ \sum _{k=1}^{p}\sum _{j=2}^{M}\Vert
(\mathfrak{L}_{i,j}^{k})^{a}u_{i,j}^{k}(\nu _{k},\phi _{k})\Vert
_{0,\widehat{\Omega }_{i,j}^{k}}^{2}\right.\nonumber\\[.2pc]
&\quad\ +\sum _{k=1}^{p}\sum _{\gamma _{s}\subseteq \Omega
^{k}}(\Vert [u]\Vert _{0,\widehat{\gamma
}_{s}}^{2}+\Vert [u_{\nu _{k}}]\Vert
_{1/2,\widehat{\gamma }_{s}}^{2}+\Vert [u_{\phi
_{k}}]\Vert _{1/2,\widehat{\gamma }_{s}}^{2})\nonumber\\[.2pc]
&\quad\ +\sum _{l\in \mathcal{D}}\sum _{k=l-1}^{l}\sum _{\gamma
_{s}\subseteq \partial \Omega ^{k}\bigcap \Gamma _{l},\mu
(\widehat{\gamma }_{s})<\infty }(\Vert u\Vert
_{0,\widehat{\gamma }_{s}}^{2}+\Vert u_{\nu _{k}}\Vert
_{1/2,\widehat{\gamma }_{s}}^{2})\nonumber\\[.2pc]
&\quad\ +\sum _{l\in \mathcal{N}}\sum _{k=l-1}^{l}\sum _{\gamma
_{s}\subseteq \partial \Omega ^{k}\bigcap \Gamma _{l},\mu
(\widehat{\gamma }_{s})<\infty }\left\Vert \left(\frac{\partial u}{\partial
n}\right)_{\widetilde{A}^{k}}^{a}\right\Vert _{1/2,\widehat{\gamma
}_{s}}^{2}\nonumber \\[.2pc]
&\quad\ +\sum _{k=1}^{p}\sum _{\gamma _{s}\subseteq B_{\rho }^{k}}(\Vert
[u]\Vert _{0,\widehat{\gamma }_{s}}^{2}+\Vert [(u)_{\nu _{k}}^{a}]\Vert
_{1/2,\widehat{\gamma }_{s}}^{2}+\Vert [(u)_{\phi _{k}}^{a}]\Vert
_{1/2,\widehat{\gamma }_{s}}^{2})\nonumber\\[.2pc]
&\quad\ +\sum _{l=1}^{L}\Vert
(\mathfrak{L}_{l}^{p+1})^{a}u_{l}^{p+1}(\xi ,\eta
)\Vert _{0,S}^{2}\nonumber\\
&\quad\ +\sum _{\gamma _{s}\subseteq \Omega ^{p+1}}(\Vert [u]\Vert
_{0,\gamma _{s}}^{2}+\Vert [u_{x_{1}}]^{a}\Vert _{1/2,\gamma
_{s}}^{2}+\Vert [u_{x_{2}}]^{a}\Vert _{1/2,\gamma
_{s}}^{2})\nonumber\\[.2pc]
&\quad\ +\sum _{l\in \mathcal{D}}\sum _{\gamma _{s}\subseteq \partial
\Omega ^{p+1}\bigcap \Gamma _{l}}\left(\Vert u\Vert _{0,\gamma
_{s}}^{2}+\left\Vert \left(\frac{\partial u}{\partial T}\right)^{a}\right\Vert _{1/2,\gamma
_{s}}^{2}\right)\nonumber\\[.2pc]
&\quad\ + \sum _{l\in \mathcal{N}}\sum _{\gamma _{s}\subseteq
\partial \Omega ^{p+1}\bigcap \Gamma _{l}}\left\Vert \left(\frac{\partial
u}{\partial N}\right)_{A}^{a}\right\Vert _{1/2,\gamma
_{s}}^{2} \left. \begin{array}{@{}c@{}}\\ \\ \\[-.3pc] \end{array} \right\}\label{eq:E-fin-quad-form}
\end{align}
We can now state the main result of this section.

\begin{theor}[\!]\label{4STN2allu} Let $\mathcal{V}^{M,W} (\{
u_{i,j}^{k}(\nu _{k},\phi _{k})\} _{i,j,k},\{
u_{l}^{p+1}(\xi ,\eta )\} _{l})$ be as defined in
{\rm (\ref{eq:E-fin-quad-form})}. Then for $M$ and $W$ large enough the
estimate
\begin{align}
&\sum_{k=1}^{p}\sum _{i=1}^{I_{k}}\left(|u_{i,1}^{k}|^{2}+\sum
_{j=2}^{M}\Vert u_{i,j}^{k}(\nu _{k},\phi _{k})\Vert _{2,\widehat{\Omega
}_{i,j}^{k}}^{2}\right)+\sum _{l=1}^{L}\Vert u_{l}^{p+1}(\xi ,\eta
)\Vert _{2,S}^{2}\nonumber\\[.2pc]
&\leq C_{M,W}\mathcal{V}^{M,W}(\{ u_{i,j}^{k}(\nu _{k},\phi _{k})\}
_{i,j,k},\{ u_{l}^{p+1}(\xi ,\eta )\} _{l})\label{eq:E-gen-stab-thm}
\end{align}
holds for all $\{ \{ u_{i,j}^{k}(\nu _{k},\phi _{k})\}_{i,j,k},\{
u_{l}^{p+1}(\xi ,\eta )\} _{l}\} \in \pi ^{M,W}$.

Here $C_{M,W}=C$ maximum $(M^{4},(\ln W)^{2})$
if there is a vertex $A_{j}$ such that Neumann boundary conditions
are imposed on the adjoining sides $\Gamma _{j}$ and $\Gamma _{j+1}$
and $C_{M,W}=C(\ln W)^{2}$ otherwise. $C$ is a constant,
independent of $M$ and $W$.
\end{theor}

Adding a weighted combination of (\ref{EabsN012upp1l}),
(\ref{EabsN2unufi}) and (\ref{EN1allurhsi}) and using the techniques and
results of \cite{pdstrk1} the result follows.\hfill $\cd$

\begin{rem} {\rm The stability theorem~3.2 holds provided the coefficients
of the differential operator $\in C^{3}(\overline{{\Omega }})$
and the curves $\phi _{i,j,l}^{k}$ and $\psi _{i,j,l}^{k}$ defined
in (\ref{eqn2.6a}), (\ref{eqn2.6b}) satisfy (\ref{eqn2.7}) for $t=1,\ldots,3$.}
\end{rem}

For problems with mixed boundary conditions the factor multiplying
the right-hand side of (\ref{eq:E-gen-stab-thm}) grows rapidly with
$M$. This creates difficulties in parallelizing the numerical scheme.
To overcome this we make the spectral element functions continuous
at the vertices of the elements. Let $\pi _{V}^{M,W}$ denote the
space of spectral element functions which are continuous at the vertices
of their elements. We define $\pi _{0}^{M,W}$ to be the space of
spectral element functions which vanish at the vertices of their element.
We now need to state a version of Theorem \ref{4STN2allu} when the
spectral element functions vanish at the vertices of their elements.

To do so, we have to prove the following result.

\setcounter{theore}{4}
\begin{lem}\label{L-sef0} Let $u_{i,j}^{k}(\xi ,\eta )$ be a
polynomial of degree $W$ in $\xi $ and $\eta $ separately{\rm ,} defined on
the unit square $S=(0,1)\times (0,1)${\rm ,} and which
is zero at all the vertices of the square. Then there exists a positive
constant $C$ such that
\begin{equation}
|u_{i,j}^{k}(\xi ,\eta )|_{0,S}^{2}\leq C(|u_{i,j}^{k}(\xi ,\eta
)|_{1,S}^{2}+|u_{i,j}^{k}(\xi ,\eta )|_{2,S}^{2}).\label{eq:E-lem-sef0}
\end{equation}
\end{lem}

Consider $u_{i,j}^{k}(\xi ,\eta )$ defined on $(0,1)\times (0,1)$. Now
$u_{i,j}^{k}(0,0)=0$. Hence
\begin{equation*}
u_{i,j}^{k}(\xi ,0)=\int _{0}^{\xi }\frac{\partial
u_{i,j}^{k}}{\partial \xi '}(\xi ',0){\rm d}\xi '.
\end{equation*}
And so we can conclude that
\begin{equation*}
|u_{i,j}^{k}(\xi ,0)|^{2}\leq \xi \int
_{0}^{1}\left|\frac{\partial u_{i,j}^{k}}{\partial \xi }(\xi
,0)\right|^{2}{\rm d}\xi.
\end{equation*}
Integrating the above with respect to $\xi $ we obtain
\begin{align}
\int _{0}^{1}|u_{i,j}^{k}(\xi ,0)|^{2}{\rm d}\xi &\leq
\frac{1}{2}\int _{0}^{1}\left|\frac{\partial u_{i,j}^{k}}{\partial \xi
}(\xi ,0)\right|^{2}{\rm d}\xi\nonumber\\[.2pc]
&\leq K(|u_{i,j}^{k}(\xi ,\eta
)|_{1,S}^{2}+|u_{i,j}^{k}|_{2,S}^{2})
\label{eq:ukij-xi0-bdd}
\end{align}
by the trace theorem for Sobolev spaces. Again
\begin{equation*}
u_{i,j}^{k}(\xi ,\eta )=u_{i,j}^{k}(\xi ,0)+\int
_{0}^{\eta }\frac{\partial u_{i,j}^{k}}{\partial \eta '}(\xi ,\eta
'){\rm d}\eta '.
\end{equation*}
Therefore
\begin{equation*}
|u_{i,j}^{k}(\xi ,\eta )|^{2}\leq
2|u_{i,j}^{k}(\xi ,0)|^{2}+2\eta \int
_{0}^{1}\left|\frac{\partial u_{i,j}^{k}}{\partial \eta }(\xi ,\eta
)\right|^{2}{\rm d}\eta.
\end{equation*}
Integrating the above with respect to $\xi$ and $\eta$ we get
\begin{align*}
\int _{S}\int |u_{i,j}^{k}(\xi ,\eta )|^{2}{\rm
d}\xi {\rm d}\eta &\leq 2\int _{0}^{1}|u_{i,j}^{k}(\xi
,0)|^{2}{\rm d}\xi\\[.2pc]
&\quad\ +\int _{S}\int \left|\frac{\partial
u_{i,j}^{k}}{\partial \eta }(\xi,\eta )\right|^{2}{\rm d}\xi
{\rm d}\eta.
\end{align*}
Combining the above with (\ref{eq:ukij-xi0-bdd}) we obtain the required
result.\hfill $\cd$\vspace{.5pc}

Clearly Lemma~\ref{L-sef0} applies equally well to any of the function
elements $u_{i,j}^{k}(\nu _{k},\phi _{k})$ for $2\leq j\leq
M$, $1\leq i\leq I_{k}$, $1\leq k\leq p$, although with a constant
$C_{k}$ which depends on $k$. Taking the supremum over the constant
$C_{k}$ (as given in (\ref{eq:E-lem-sef0})) we conclude that
\begin{equation}
|u_{i,j}^{k}(\nu _{k},\phi _{k})|_{0,\widehat{\Omega
}_{i,j}^{k}}^{2}\leq C(|u_{i,j}^{k}(\nu _{k},\phi
_{k})|_{1,\widehat{\Omega }_{i,j}^{k}}^{2}+|u_{i,j}^{k}(\nu _{k},\phi
_{k})|_{2,\widehat{\Omega }_{i,j}^{k}}^{2}),\label{eq:ukij-nuk-phik-bdd}
\end{equation}
for all function elements with $1\leq k\leq p, 1\leq i\leq I_{k}$,
$2\leq j\leq M$. Here $C$, of course, denotes a generic constant. We can
now state the final result of this section.

\setcounter{theore}{2}
\begin{theor}[\!]\label{T-stab-est-v0} Let $\{ \{ u_{i,j}^{k}(\nu
_{k},\phi _{k})\} _{i,j,k},\{ u_{l}^{p+1}(\xi ,\eta )\} _{l}\} $ belong
to the space of functions $\pi _{0}^{M,W}$ which are zero at the
vertices of the elements on which they are defined. Then the following
estimate holds{\rm :}
\begin{align}
&\sum_{k=1}^{p}\sum _{j=2}^{M}\sum _{i=1}^{I_{k}}\Vert u_{i,j}^{k}(\nu
_{k},\phi _{k})\Vert _{2,\widehat{\Omega }_{i,j}^{k}}^{2}+\Vert
u_{l}^{p+1}(\xi ,\eta )\Vert _{2,S}^{2}\nonumber\\[.2pc]
&\leq C(\ln W)^{2}\mathcal{V}^{M,W}(\{ u_{i,j}^{k}(\nu _{k},\phi _{k})\}
_{i,j,k},\{ u_{l}^{p+1}(\xi ,\eta )\} _{l}) \label{eq:E-stab-est-v0}
\end{align}
for $M$ and $W$ large enough.
\end{theor}

In the above $u_{i,1}^{k}(\nu _{k},\phi _{k})$ is taken to be
identically zero for $1\leq k\leq p$ and $1\leq i\leq I_{k}$.

\noindent Combining the estimates (\ref{eq:E-lem-sef0}) and
(\ref{eq:ukij-nuk-phik-bdd}) with the earlier results
(\ref{eq:E-stab-est-v0}) follows.\hfill $\cd$

\section{\label{sec:Gen-Conc}Conclusion}

We can use the stability theorem \ref{4STN2allu} to formulate a
numerical scheme to obtain an approximate solution to the elliptic
boundary value problem (\ref{eq:Eelloper}) as has been described in
\cite{tomar-01,tomar-dutt-kumar-02}. For problems with Dirichlet
boundary conditions we choose our solution to be a non-conforming
spectral element representation which minimizes a functional which is
the sum of the squares of weighted squared norms of the residuals in the
partial differential equation and fractional Sobolev norms of the
residuals in the boundary conditions and a term which measures the sum
of the jumps in the function and its derivatives in appropriate Sobolev
norms at inter-element boundaries. In a sectoral neighbourhood of the
corners these quantities are computed using modified polar coordinates
and in the remaining part of the domain we use a global coordinate
system. This method is faster than the h-p finite element method
as there are no common boundary values to solve for
\cite{tomar-01,tomar-dutt-kumar-02}.

For problems with mixed boundary conditions we have to make the spectral
element functions continuous only at the vertices of the elements. As a
result the Schur complement matrix has a small dimension and an accurate
inverse can be computed. Hence the numerical scheme has a computational
complexity which is less for finite element methods.

Moreover, the construction of a pre-conditioner for the Schur complement
matrix is very simple unlike the case for finite element methods. In
fact, for problems in three dimensions the construction of
pre-conditioners for the Schur complement matrix becomes quite complex
for finite element methods \cite{pavwid}.

Though the ideas in these papers deal with problems in two dimensions,
they generalize to three dimensions. We intend to study these problems
both theoretically and computationally in future work.

\section*{Acknowledgements}

The research by PD is partly supported by Aeronautical Research and
Development Board (ARDB) and Center for Development of Advanced
Computing (CDAC), Pune, India. The research by ST is supported by
Council for Scientific and Industrial Research (CSIR), India.

\end{document}